\documentclass[hidelinks,onefignum,onetabnum]{siamart250106}

\usepackage{lipsum}
\usepackage{amsfonts}
\usepackage{graphicx}
\usepackage{epstopdf}
\usepackage{algorithmic}

\usepackage{amsmath,amssymb,mathtools,upgreek}
\usepackage{xcolor,soul,bm}

\newcommand{\red}[1]{\textcolor{red}{#1}}

\usepackage[shortlabels]{enumitem}

\usepackage{booktabs,multirow,geometry,array,makecell,xfrac}
\usepackage[flushleft]{threeparttable}

\DeclareMathOperator{\tr}{tr}
\DeclareMathOperator{\Ad}{Ad}
\DeclareMathOperator{\ad}{ad}
\DeclareMathOperator{\Cov}{\mf\Sigma}
\DeclareMathOperator{\vecn}{vec}

\newcommand{\mf}[1]{\mathbf{#1}}
\newcommand{\mc}[1]{\mathcal{#1}}
\usepackage[scr=boondoxo]{mathalpha}
\newcommand{\ms}[1]{\mathscr{#1}}

\newcommand*{\tran}{^{\mkern-1.5mu\mathsf{T}}}
\newcommand*{\htran}{^{\mkern-0.5mu\mathsf{H}}}

\newcommand{\ip}[2]{\langle{#1},{#2}\rangle}

\newcommand{\rexp}[1]{\textup{Exp}_{#1}}
\newcommand{\rlog}[1]{\textup{Log}_{#1}}

\newcommand{\bmu}{\bm\upmu}

\newcommand{\subs}[1]{\textup S_{#1}}

\newcommand{\ii}{\mathrm{i}}
\newcommand{\LE}{\scriptscriptstyle{LE}}
\newcommand{\conv}{\ast}

\renewcommand{\proofbox}{\rule{1.5ex}{1.5ex}}

\ifpdf
  \DeclareGraphicsExtensions{.eps,.pdf,.png,.jpg}
\else
  \DeclareGraphicsExtensions{.eps}
\fi

\newsiamremark{remark}{Remark}
\newsiamremark{hypothesis}{Hypothesis}
\crefname{hypothesis}{Hypothesis}{Hypotheses}
\newsiamthm{claim}{Claim}

\headers{Means of Random Variables in Lie Groups}{Shiraz Khan, Jikai Ye, and Gregory S. Chirikjian}

\title{Means of Random Variables in Lie Groups
\thanks{Submitted on July 11, 2025
}}

\author{Shiraz Khan\thanks{Department of Mechanical Engineering, University of Delaware, Newark, DE 
  (\email{mail@shiraz-k.com}).} \and 
  Jikai Ye\thanks{Department of Mechanical Engineering, National University of Singapore, Singapore 
  (\email{e0925572@u.nus.edu}).}
\and Gregory S. Chirikjian\thanks{Department of Mechanical Engineering, University of Delaware, Newark, DE 
  (\email{gchirik@udel.edu}).}}

\usepackage{amsopn}

\ifpdf
\hypersetup{
  pdftitle={Means of Random Variables in Lie Groups},
  pdfauthor={S. Khan, J. Ye, and G. S. Chirikjian}
}
\fi

\begin{document}

\maketitle

\begin{abstract}
The concepts of mean (i.e., average) and covariance of a random variable are fundamental in statistics, and are used to solve real-world problems such as those that arise in robotics, computer vision, and medical imaging. 
On matrix Lie groups, multiple competing definitions of the mean arise, including the Euclidean, projected, distance-based (i.e., Fréchet and Karcher), group-theoretic, and parametric means. This article provides a comprehensive review of these definitions, investigates their relationships to each other, and determines the conditions under which the group-theoretic means minimize a least-squares type cost function. We also highlight the dependence of these definitions on the choice of inner product on the Lie algebra.
The goal of this article is to guide practitioners in selecting an appropriate notion of the mean in applications involving matrix Lie groups.
\end{abstract}



\section{Introduction}
\subsection{Background and Motivation}
The concept of the \textit{mean} (i.e., arithmetic average) of a random vector has been used for statistical applications since the works of Legendre and Gauss, who observed that the mean of a set of vectors minimizes the sum of its squared Euclidean distances to the vectors. 
A modern and more sophisticated notion of the mean was studied by Élie Cartan \cite{barbaresco2012information,arnaudonMediansMeansRiemannian2013a}, who was interested in defining the mean of a set of samples in a Riemannian manifold. As there is generally no canonical choice of Riemannian metric for a given manifold, there is no canonical definition of the mean on an arbitrary manifold. Hence, for statistical applications involving non-Euclidean manifolds, one must choose between the various definitions of means based on additional considerations. One such consideration is the existence of a group structure on the manifold, which leads to the notion of a Lie group \cite{chirikjian2011stochastic,barfoot2024state,gallier2020differential}.
The Lie groups that are most commonly encountered in applications are \textit{matrix Lie groups},
and are represented using a set of invertible matrices that is closed under the multiplication and inversion of matrices (which are analytic maps). 
The concept of the mean in matrix Lie groups is used in state estimation \cite{long2013banana,chahbazian2022generalized,brossard2017unscented}, generative modeling \cite{jagvaral2024unified,jiang2023se}, and imaging \cite{xu2020lightweight}, among other applications.
One expects the mean of a random variable in a matrix Lie group to transform in a natural way (as will be made precise in this article) when the underlying random variable is transformed. This ensures that when the choice of units or coordinate system is changed, the mean of a random variable transforms in accordance. 
There are other considerations that are taken into account when choosing a definition of the mean for a given application, such as existence and uniqueness, ease of computation, and the ability to quantify uncertainty (e.g., via a generalized notion of the covariance matrix).

\subsection{Literature Review}
In the past, several works have investigated the problem of defining the \textit{mean} of a random variable that assumes values in a non-Euclidean space, such as a Riemannian manifold or a Lie group \cite{nielsen2013matrix}. 
Means on homogeneous spaces such as spheres and symmetric positive definite (SPD) matrices have also been studied \cite{barbaresco2012information,buss2001spherical,moakher2005differential}. However, the aforementioned works restrict their discussion to the case where a discrete set of samples are obtained from a random variable, failing to discuss the extension to continuous random variables on Lie groups, and the ensuing details related to the Haar measure of the group. 
Continuous random variables on Lie groups can be represented using parametric probability density functions (pdfs) \cite{long2013banana,said2017gaussian,jupp1989unified} or by using Monte Carlo methods such as particle filtering \cite{zhang2017feedback}.
With the advent of deep learning tools like normalizing flows, Moser flows, and diffusion models, it is now possible to represent fairly complicated pdfs on Lie groups
\cite{falorsi2019reparameterizing,de2022riemannian}. Given a continuous random variable on a Lie group, one is often interested in summarizing the random variable using a mean and covariance, e.g., for downstream control and decision-making tasks \cite{long2013banana,zhang2017feedback,bordin2022diff,barfoot2014associating}. Hence, the problem of defining and computing the mean of a continuous random variable on a Lie group is of practical relevance.


A secondary limitation of the existing works is that they do not explain how the various definitions of means relate to one another. For instance, the usual Euclidean notion of the mean as well as the mean computed using the \textit{exponential parametrization} of the group (which we call the $\log$-Euclidean mean) are both commonly-used by practitioners \cite{ye2024uncertaintyCDC,bourmaud2015continuous,ge2024geometric}. However, the literature that discusses the intrinsic (e.g., Riemannian and group-theoretic) notions of the mean often omits the extrinsic (e.g., Euclidean and parametric) means in its discussion, as the latter are coordinate-dependent. Nevertheless, a clear and careful exposition of the differences between the various competing
definitions of means -- without extensively relying on differential geometric concepts -- would benefit practitioners in disciplines like robotics and computer science who are already acquainted with the language of Lie groups.

A third point which this article aims to address is that of the dependence of the definitions of means on the choice of basis and/or inner product on the Lie algebra. It is well-known that an $\Ad$-invariant inner product on the Lie algebra (if one exists) makes the group-theoretic definition of the mean coincide with the Riemannian ones. However, there is little discussion of what happens when one works with an inner product that is not $\Ad$-invariant. For instance, the geodesics corresponding to an $\Ad$-invariant inner product on $\mathfrak{so}(3)$ (i.e., the Lie algebra of $SO(3)$) are rotations about a fixed axis, but the same is not true when one works with a weighted inner product on $\mathfrak {so}(3)$. 
There is a similar gap in the literature when it comes to the mean of a random variable in $SE(3)$ (i.e., the Special Euclidean group, which comprises all the rigid transformation of $\mathbb R^3$); while the works in \cite{zefran1999metrics,miolane2015computing,diatta2024dual} have studied the Riemannian exponential map of $SE(3)$, the implications of their findings to the problem of defining a mean for random variable in $SE(3)$ have not been explored.

\subsection{Contributions and Overview}
In this article, we review, compare, and contrast the existing definitions of means and covariances for a random variable that takes values in a matrix Lie group. The dependence of each of these definitions on additional structure (e.g., a choice of inner product and/or basis on the Lie algebra) is clarified. In particular, \Cref{thm:lie-karcher,thm:symmetric} address the question of whether the group-theoretic definition of mean (which is defined using the Lie-theoretic log map) minimizes a cost function similar to the Euclidean, Fr\'echet, and Karcher means. The examples of $SO(d)$ and $SE(d)$ (where $d=2$ or $3$) are used to illustrate the main ideas. Finally, the appendices explain the technical aspects of defining and working with pdfs on Lie groups, keeping the discussion accessible to practitioners.

The organization of this article is as follows. \Cref{sec:means} introduces the various definitions of means on Lie groups.
In \Cref{sec:relationship}, we clarify the relationship between the group-theoretic and Karcher means, which coincide under certain (sometimes unsatisfiable) conditions. \Cref{sec:properties} discusses the various properties that one desires of a notion of mean on a Lie group, such as its compatibility with the algebraic structure of the Lie group.
Finally, \Cref{sec:examples} discusses the examples of $SO(d)$ and $SE(d)$ in detail. To keep the main text accessible and easy to follow, the technical details of calculus and Riemannian geometry on Lie groups are deferred to the appendices.


\section{Definitions of Means}\label{sec:means}
\subsection{Euclidean Mean}
\label{sec:euclidean}
Given a set of $N$ vectors $\lbrace\mf x_i\rbrace_{i=1}^N$ in $\mathbb R^n$, 
its {Euclidean mean} is defined as the vector $\frac{1}{N}\sum_{i=1}^N \mf x_i$, which is also known as the arithmetic average.
Let $\lVert \mf x \rVert_{\mf W}\coloneqq \sqrt{\mf x\tran \mf W \mf x}$ denote the weighted Euclidean norm, where $\mf W$ is any symmetric positive definite (SPD) matrix.
It is known that the Euclidean mean has the following characterization:
\begin{align}
   \bmu_E\big(\lbrace\mf x_i\rbrace_{i=1}^N\big) \coloneqq \frac{1}{N}\sum_{i=1}^N \mf x_i= \underset{\mf y\in\mathbb R^n}{\arg\min} \sum_{i=1}^N \lVert \mf x_i - \mf y \rVert^2_{\mf W} 
    \label{eq:discrete-euclidean}
\end{align}
in which the second equality holds irrespective of the choice of (positive definite) weighting matrix, $\mf W$. 
That the second equality of \cref{eq:discrete-euclidean}
holds irrespective of the choice of weighting matrix $\mf W$ is a property that, as we will show in \Cref{sec:relationship}, does not extend to Lie groups.
Suppose we have a random variable $\tilde{\mf x}$ in $\mathbb R^n$ whose probability density function (pdf) is $f$, then we define the {Euclidean mean} of $\tilde{\mf x}$ via an integral:
\begin{align}
\bmu_E(\tilde{\mf x}) \coloneqq \int_{\mathbb R^n}\mf x f(\mf x) d\mf x =\underset{\mf y\in\mathbb R^n}{\arg\min}\ \int_{\mathbb R^n} \lVert \mf x - \mf y \rVert^2_{\mf W} f(\mf x) d\mf x.
\label{def:euclidean-mean}
\end{align}
Observe that (\ref{eq:discrete-euclidean}) arises as a special case of (\ref{def:euclidean-mean}) for the case where $f$ is given by a sum of Dirac delta functions: $f(\mf x)=\frac{1}{N}\sum_{i=1}^N \delta(\mf x - \mf x_i)$.\footnote{This can be made rigorous by considering instead of $f(\mf x) d\mf x$, a {probability measure} on $\mathbb R^n$. By replacing $f(\mf x) d\mf x$ with a \textit{counting measure} on $\lbrace\mf x_i\rbrace_{i=1}^N$, (\ref{def:euclidean-mean}) becomes  (\ref{eq:discrete-euclidean}) \cite[Sec. 5.5.1]{pennec2020beyond}.
Throughout this article, we use Dirac delta functions as a notational shorthand for converting between probability measures that are absolutely continuous w.r.t.~the integration measure (whose pdfs are continuous functions on $G$) and counting measures (whose pdf is technically not a function, but a \textit{distribution}). See Appendix \ref{app:mlg-calc} for the technical details.
\label{foot:dirac-delta}
}
The Euclidean mean is alternatively referred to as the \textit{average}, \textit{centroid}, \textit{center of mass}, or \textit{barycenter}.

It is clear that the Euclidean mean generalizes to vector spaces other than $\mathbb R^n$. For instance, let $\mathbb R^{m\times m}$ denote the space on $m\times m$ real matrices, where $m$ is a positive integer. 
If $\tilde{\mf A}$ is a matrix-valued random variable in $\mathbb R^{m\times m}$, then its Euclidean mean can be defined as
\begin{align}\label{def:euclidean-matrix-mean}
\bmu_E(\tilde{\mf A}) &\coloneqq 
\int_{\mathbb R^{m\times m}}\mf A f(\mf A) d\mf A = \int_\mathbb R\int_\mathbb R\cdots \int_\mathbb R \mf A f(\mf A)\, d A_{11}\,d A_{12}\cdots d A_{mm}, 
\end{align}
where $d\mf A$ denotes the Lebesgue measure of $\mathbb R^{m\times m}$, and $A_{ij}$ denotes the $(i,j)^{th}$ component of $\mf A$.
To see the relationship between \cref{def:euclidean-matrix-mean} and \cref{def:euclidean-mean}, consider the \textit{vectorization} map 
\begin{align*}
\vecn: \mathbb R^{m\times m} \rightarrow \mathbb R^{m^2}
\end{align*}
that stacks the components of a matrix to create a column vector. Since $\vecn(\tilde{\mf A})$ is a vector-valued random variable, we can use the linearity of $\vecn$ to show that $\bmu_E(\tilde{\mf A})=\vecn^{-1}\big(\bmu_E\big(\vecn(\tilde{\mf A})\big)\big).$
Furthermore, we have $\lVert \mf A\rVert \coloneqq 
\sqrt{\tr\big(\mf A\tran {\mf A}\big)} = \lVert \vecn({\mf A}) \rVert,$
which is the \textit{Frobenius norm} of the matrix ${\mf A}\in \mathbb R^{m\times m}$. Here, 
$\lVert\mf x \rVert\coloneqq\lVert\mf x \rVert_{\mf I_n} = \mf x\tran \mf x$ denotes the non-weighted Euclidean norm and $\mf I_n$ is the $n\times n$ identity matrix.
Using \cref{def:euclidean-mean}
and the linearity of $\vecn$, it follows that
\begin{align}\label{eq:euclidean-matrix-mean-frob}
\bmu_E(\tilde {\mf A})= \underset{\mf B\in\mathbb R^{m\times m}}{\arg\min}\ \int_{\mathbb R^{m\times m}} 
\lVert \mf A - \mf B\rVert^2\,
f(\mf A) d\mf A.
\end{align}

The above can be generalized to other weighted variants of the Frobenius norm. For instance, one can consider the weighted Frobenius norm, $\lVert \mf A \rVert_{(\mf P,\mf Q)}\coloneqq \lVert\mf P \mf A \mf Q\rVert$, where $\mf P$ and $\mf Q$ are a given pair of invertible matrices in $\mathbb R^{m\times m}$.  Since 
$\mf A \mapsto \mf P \mf A \mf Q$ is an invertible linear transformation of $\mf A$, it can be represented using an invertible matrix $\mf M \in \mathbb R^{m^2 \times m^2}$, such that $\lVert \mf A \rVert_{(\mf P,\mf Q)} = \lVert \mf M \vecn(\mf A) \rVert=\lVert \vecn(\mf A) \rVert_{\mf M\tran \mf M }$. From (\ref{def:euclidean-mean}), it follows that the matrix $\bmu_E(\tilde {\mf A})$ also minimizes the (expected value of the) weighted Frobenius norm.

Finally, we point out that the discussion above readily extends to the case where $\tilde {\mf A}$ is a random variable in $\mathbb C^{m\times m}$ -- the space of $m\times m$ matrices with complex entries. In the complex case, one defines the vectorization map $\mf A \mapsto \vecn(\mf A)$ such that it stacks the real and imaginary components of $\mf A$ into a single column, consisting of $2m^2$ real numbers. 

\subsection{Extrinsic Means}
\label{sec:extrinsic}
In this article, we are interested in the case where the random variable takes values in a matrix Lie group. Let $G$ be a real $n$-dimensional Lie group endowed with a real $m$-dimensional matrix representation, where $n$ and $m$ are positive integers.\footnote{As explained in \Cref{app:mlg}, a real Lie group may have a complex representation, such as in the case of the {Unitary group}, $U(1)$. Here, the word \textit{real} indicates that the Lie algebra of $U(1)$ is a real vector space, i.e., is closed under multiplication by real (but not complex) numbers.} The elements of $G$ constitute a subset of $\mathbb R^{m\times m}$ that is closed under matrix multiplication and inversion. 
The group multiplication $(\mf g, \mf h)\mapsto \mf g \mf h$ and inversion $\mf g \mapsto \mf g^{-1}$ are analytic functions \cite{hall2013lie}. The Lie algebra of $G$ is denoted as $\mathfrak g$, and is a subspace of $\mathbb R^{m\times m}$ that is closed under the \textit{Lie bracket} operation, defined as $[\mf X,\mf Y]\coloneqq \mf X \mf Y - \mf Y \mf X$ for all $\mf X,\mf Y\in \mathfrak g$. Note that $\mathfrak g$ is a vector space consisting of $m\times m$ matrices, such that any basis for this vector space has exactly $n$ matrices.

The focus of this article is on matrix Lie groups (reviewed in \Cref{app:mlg}) that are frequently encountered in scientific and engineering applications, such as the \textit{Special Orthogonal} group $SO(d)$ and the \textit{Special Euclidean} group $SE(d)$ (where $d=2$ or $3$). In the case of $G=SO(d)$, each element $\mf R\in SO(d)$ is a $d\times d$ orthogonal matrix (i.e., $\mf R^{-1}=\mf R\tran$), so that the dimension of the representation (i.e., the number $m$) is equal to $d$. The Lie algebra $\mathfrak{so}(d)$ consists of $d\times d$ skew-symmetric matrices, and is a vector subspace of $\mathbb R^{d\times d}$ whose dimension is $n=\frac{d(d-1)}{2}$ \cite{hall2013lie}. 
In the case of $G=SE(d)$, the elements of the group are typically represented using \textit{homogeneous transformation} matrices. These are matrices in $\mathbb R^{(d+1)\times (d+1)}$ of the form
\begin{align}
\mf H(\mf R, \mf t) \coloneqq
\left[
\begin{array}{@{}lc@{}}
\,\mf R & \mf t \\
\,\mf 0_{1\times d} & 1
\end{array}
\right] \in SE(d),
\end{align}
where $\mf R\in SO(d)$, $\mf t\in\mathbb R^d$, and $\mf 0_{1\times d}$ is the $1\times d$ matrix of zeros. The multiplication of homogeneous transformations can be evaluated using matrix multiplication, as follows:
\begin{align}
    \mf H(\mf R,\mf t)\,\mf H(\mf Q,\mf s)=\mf H(\mf R\mf Q,\mf R\mf s + \mf t),
    \label{eq:se3-multiplication}
\end{align}
where $\mf R,\mf Q\in SO(d)$, and $\mf t, \mf  s\in\mathbb R^d$.
The dimension of $\mathfrak{se}(d)$ is $n=\frac{d(d+1)}{2}$. Thirdly, we let $GL(d,\mathbb R)$ denote the matrix Lie group consisting of $d\times d$ invertible matrices with real entries, called the \textit{General Linear} matrix Lie group, whose dimension is $n=d^2$. Finally, $(\mathbb R^n,+)$ will denote the space $\mathbb R^n$ equipped with vector addition as the group operation, which is an Abelian (i.e., commutative) Lie group.\footnote{While $(\mathbb R^n,+)$ can also be represented using matrices (e.g., as matrices of the form $\mf H(\mf I_n, \mf x)$, where $\mf x\in \mathbb R^n$), we will denote the elements of $(\mathbb R^n,+)$ as vectors for clarity. For instance, the inverse of the element $\mf x\in(\mathbb R^n,+)$ will not be written as ``$\mf x^{-1}$", but as $-\mf x$. Note that the Lie algebra of $(\mathbb R^n, +)$ is also $\mathbb R^n$, and the Lie bracket vanishes identically.}

Consider a random variable $\tilde {\mf g}$ that takes values in $G$; we will denote this situation as $\tilde {\mf g}\in G$. 
The pdf of $\tilde {\mf g}$, denoted as $f$, may be defined w.r.t.~the left Haar measure of $G$ (see \Cref{app:mlg-calc} for a review). We write `$d\mf g$' to denote integration with the left Haar measure, so that $\int_G f(\mf g) d\mf g = 1$. 
With these definitions, the \textit{extrinsic Euclidean mean} of $\tilde {\mf g}$ can be defined as follows:
\begin{align}\label{def:group-euclidean}
    \bmu_E\big(\tilde{\mf g}\big) &=\int_{G}  \mf g\,
f(\mf g)\,d\mf g .
\end{align}
Given an arbitrary (fixed, deterministic) vector $\mf x \in \mathbb R^m$, the quantity $\tilde{\mf g}\,\mf x$ can be viewed as a random vector in $\mathbb R^m$. One can use the linearity of the Haar integral to show that
    $
    \bmu_E\big(\tilde{\mf g}\,\mf x\big)= \bmu_E(\tilde{\mf g}) \mf x,
    $
    which gives an alternative interpretation for $\bmu_E(\tilde{\mf g})$; it describes the Euclidean mean of the vector-valued random variable, $\tilde{\mf g}\,\mf x$. 
    

The main limitation of the Euclidean mean is that $\bmu_E(\tilde{\mf g}) \notin G$ in general.
For instance, the Euclidean mean of a random variable on $SO(2)$ will generally not be a rotation matrix.
To ensure that the mean is an element of $G$, one can construct a function $
\mc P:\mathbb R^{m\times m} \rightarrow G$ that {projects} vectors onto $G$, and define the \textit{projected mean} of $\tilde {\mf g}$ as $\mc P\big(\bmu_E(\tilde{\mf g}) \big)$ \cite{fletcher2003statistics}.
In the case of $\tilde{\mf g} \in SO(2)$, we can consider the projection function\footnote{Note that this function is ill-defined on a measure-zero subset of $\mathbb R^{2\times 2}$ where $\mf A$ fails to be invertible. Nevertheless, this issue is easily resolved in practice by adding an arbitrarily small amount of noise to the data.\label{foot:projection}} $\mc P_{SO}(\mf A) \coloneqq (\mf A \mf A\tran)^{-\frac{1}{2}}\mf A=\mf A(\mf A\tran \mf A)^{-\frac{1}{2}}$, which picks out the orthogonal matrix in the polar decomposition of $\mf A$ \cite[Thm. 7.3.1]{hornMatrixAnalysis2012}.
Most of the extrinsic definitions of \textit{mean} that are encountered in the literature can be viewed as a projected mean for a suitable choice of projection function $\mc P$ \cite{fletcher2003statistics,srivastava2002monte}.

\subsection{Distance-Based Means}
\label{sec:distance}
The Euclidean mean of a random vector, defined via (\ref{def:euclidean-mean}), can be generalized as follows. Let $G$ be equipped with a {distance function}
$\mc D: G \times G \rightarrow \mathbb R_{\geq 0}$.\footnote{ 
We say that $\mc D$ is a \textit{distance function} if it satisfies, for all $\mf g_1,\mf g_2, \mf g_3\in G$, the following four properties (i) $\mc D(\mf g_1, \mf g_1) = 0$, (ii) $\mc D(\mf g_1, \mf g_2) > 0$ when $\mf g_1 \neq \mf g_2$, (iii) symmetry: $\mc D(\mf g_1, \mf g_2)=\mc D(\mf g_2, \mf g_1)$, and (iv) triangle inequality:
$\mc D(\mf g_1, \mf g_2) + \mc D(\mf g_2, \mf g_3)\geq \mc D(\mf g_1, \mf g_3).$ The term \textit{distance function} is used in place of the term \textit{metric} in order to differentiate it from the notion of a Riemannian metric (c.f. \Cref{app:riemannian}). We assume that $(G, \mc D)$ is a \textit{complete} metric space, so that the minimizer in (\ref{eq:distance-mean-functional}) is a point in $G$.
} 
Then, we can consider the minimizers of the following optimization problem:
\begin{align}
    \underset{\mf h \in G}{\textrm{minimize}}\left(
    \int_{G}\mc D(\mf g, \mf h)^2 f(\mf g)\,d \mf g\right).
    \label{eq:distance-mean-functional}
\end{align}
A minimizer of (\ref{eq:distance-mean-functional}) is called a \textit{Fréchet mean}. We let $\mu_F(\tilde {\mf g};\mc D)\subseteq G$ denote the \textit{set} of Fréchet means of $\tilde {\mf g}$, explicitly indicating its dependence on the choice of distance function $\mc D$. In the special case of $G=(\mathbb R^n,+)$ and $\mc D(\mf x, \mf y)=\lVert \mf x - \mf y\rVert_{\mf W}$ (where $\mf W$ is any SPD matrix), the Fr\'echet mean is unique and can be described using the closed-form expression in \cref{def:euclidean-mean}. If $G$ is equipped with a Riemannian metric (reviewed in Appendix \ref{app:riemannian}), then the corresponding \textit{Riemannian distance} function $\mc R$ can be used to define a Fréchet median or mean, where $\mc R(\mf g, \mf h)$ is the length of shortest (i.e., length-minimizing) geodesic connecting $\mf g$ and $\mf h$. We refer the corresponding set of Fr\'echet means as the \textit{Riemannian Fréchet means}, denoted as $\mu_{F}(\tilde {\mf g};\mc R)$.

In recent literature, the term \textit{Karcher mean} has come to refer to a point $\bmu\in G$ that {locally} minimizes (\ref{eq:distance-mean-functional}) with $\mc R$ as the choice of distance function \cite{afsariRiemannianCenter2011,arnaudonMediansMeansRiemannian2013a,pennec2013bi}. Let $\rlog{\mf h}$ denote the Riemannian log map at $\mf h$. The vector $\rlog{\mf h}(\mf g)$ is an element of $T_{\mf h} G$ (the tangent space of $G$ at $\mf h$), and represents the initial velocity of the shortest geodesic from $\mf h$ to $\mf g$.\footnote{If there does not exist a geodesic connecting $\mf g$ and $\mf h$ (e.g., if $\mf g$ and $\mf h$ lie in disconnected components of $G$), then $\mc R(\mf g, \mf h)$ is given the value $+\infty$ by convention, and the map $\rlog{\mf h}$ is undefined at $\mf g$. For this reason, $\rlog{\mf h}$ is referred to as a \textit{local} inverse of the Riemannian exponential map, $\rexp{\mf h}$  \cite{lee2018introduction}.} Note that $\mf h^{-1} \rlog{\,\mf h}(\mf g)\in \mathfrak g$.
Using the first-order condition for local minimization of (\ref{eq:distance-mean-functional}) (i.e., setting the gradient to zero) \cite[Thm. 1.2]{karcherRiemannianCenterMass1977}, we say that $\bmu \in G$ is a Karcher mean if
\begin{align}
    \int_{G'} \rlog{\bmu} (\mf g) f(\mf g) d\mf g = \mf 0_m,
    \label{eq:karcher-definition}
\end{align}
where $\mf 0_m$ is the $m\times m$ matrix of zeros and $G' \subseteq G$ is a neighborhood of $\bmu$ 
that satisfies the following conditions:
\begin{enumerate}[(i)]
    \item $\int_{G'}f(\mf g)d\mf g=1$, i.e., the support of $f$ is essentially contained in $G'$, and
    \item $\rlog{\bmu}(\mf g)$ is well-defined for all $\mf g \in G'$.
\end{enumerate}
On account of property (ii), $G'$ is called a \textit{normal neighborhood} of $\bmu$ in $G$ \cite{lee2018introduction}.
Conditions for the uniqueness and existence of a Karcher mean were studied by Karcher and Grove \cite{grove1973conjugate,karcherRiemannianCenterMass1977} and several authors have since studied generalized versions of the same concept \cite{afsariRiemannianCenter2011,arnaudonMediansMeansRiemannian2013a,pennec2012exponential}.%
\footnote{In the literature that studies the uniqueness and existence of Karcher means, a stronger condition is placed on the set $G'$, namely, that it is \textit{geodesically convex}. This requires any two points of $G'$ to be connected by a unique geodesic that is contained in $G'$. Nevertheless, our (relaxed) characterization of $G'$ is sufficient for (\ref{eq:karcher-definition}) to be well-defined, and serves to make the definition of the Karcher means similar to that of the \textit{group-theoretic means} introduced in Section \ref{sec:group-theoretic}.}
We let $\mu_K(\tilde {\mf g};\mc R)$ denote the set of all Karcher means. In \Cref{sec:relationship}, we discuss some special cases where the functions $\mc R$ and $\rlog{\mf h}$ have closed-form expressions. However, in general, these functions are defined as solutions to first order ordinary differential equations (ODEs), as explained in \Cref{app:geodesics}.

\subsection{Group-Theoretic Means}
\label{sec:group-theoretic}
The group-theoretic exponential map of $G$, denoted as $\exp$, assigns to each vector in $\mathfrak g$ an element in $G$ (see App. \ref{app:mlg} for a review). 
In general, this map is neither injective nor surjective.
Nevertheless, there always exists a subset $\subs{\mathfrak g}\subseteq \mathfrak g$ containing the origin of $\mathfrak g$ on which the exponential map $\exp:\subs{\mathfrak g} \rightarrow \subs{G}$ and its inverse $\log:\subs{G} \rightarrow \subs{\mathfrak g}$ are well-defined, with $\subs{G}=\exp(\subs{\mathfrak g})\subseteq G$ \cite{hall2013lie}. The map $\log$ is called the \textit{principal logarithm} of $G$ and is (in general) different from the Riemannian $\rlog\bmu$ map used to define the Karcher mean. Hence, the group-theoretic $\log$ map can be used to define a notion of mean on $G$ that (in general) differs from the Riemannian Fréchet and Karcher means \cite{chirikjian2000engineering,wang2008nonparametric}. Given a point $\bmu \in G$, we say that $\bmu$ is a \textit{group-theoretic mean} of $\tilde{\mf g}$ if
\begin{align}\label{def:group-theoretic-mean}
\int_{G'} \log(\bmu^{-1} \mf g) f(\mf g) d \mf g = \mf 0_m,
\end{align}
where $G' \subseteq G$ is a neighborhood of $\bmu$ that has the following properties:
\begin{enumerate}[(i)]
    \item $\int_{G'}f(\mf g)d\mf  g=1$, i.e., the support of $f$ is (essentially) contained in $G'$, and
    \item $\log(\bmu^{-1}\mf g)$ is well-defined for all $\mf g \in G'$. 
\end{enumerate}
Hereafter, we write $G'$ to denote a subset of $G$ satisfying the foregoing properties. The set of all group-theoretic means is denoted by $\mu_G(\tilde{\mf g})$.%
\footnote{
    The set $G'$ is typically chosen as $\bmu \subs{G}$, where $\bmu \subs{G}\coloneqq \lbrace \bmu \mf g\,|\,\mf g\in \subs{G}\rbrace$ is a neighborhood of $\bmu$, and $\subs{G}$ is the domain of the $\log$ map. This satisfies condition (ii), and only condition (i) remains to be checked. Condition (i) is automatically satisfied for Abelian and compact Lie groups, as well as some non-compact Lie groups like $SE(2)$ and $SE(3)$, since the set $\subs{G}$ for these groups contains all but a measure-zero subset of $G$.
    
    An example of a Lie group for which condition (i) is restrictive is $SL(2, \mathbb R)$ (c.f. \cite[Remark 1]{pennec2020beyond}), in which case it is possible for $\mu_G(\tilde{\mf g})$ to be the empty set. 
}
When $G=(\mathbb R^n, +)$, it holds that $\log(\mf y^{-1}\mf x)= \mf x - \mf y$ for all $\mf x,\mf y \in \mathbb R^n$. Given an $\mathbb R^n$-valued random variable $\tilde {\mf x}$, we have $\mu_G(\tilde {\mf x}) = \lbrace\bmu_E(\tilde {\mf x})\rbrace$, which shows that $\mu_G$ is a natural generalization of the Euclidean mean from a group-theoretic standpoint.

Let $\Ad$ and $\ad$ denote the adjoint representations of $G$ and $\mathfrak g$, respectively (reviewed in Appendix \ref{app:mlg}).
Using the identity $$\Ad_{\mf g} \big(\log(\mf h)\big) = \mf g \log(\mf h)\mf g^{-1}= \log(\mf g \mf h \mf g^{-1})\quad \forall \mf g,\mf h\in G,$$
and the linearity of $\Ad_{\bmu}$,
we can restate \cref{def:group-theoretic-mean} as
\begin{align}
\Ad_{\bmu} \left( \int_{G'} \log(\bmu^{-1} \mf g) f(\mf g) d \mf g \right) 
&=\int_{G} \log( \mf g\bmu^{-1}) f(\mf g) d \mf g=\mf 0_m.
\label{def:group-theoretic-mean-right}
\end{align}

\subsection{Parametric Means}\label{sec:log-euclidean}
Let $\subs{\mathbb R^n}$ be an open subset (viewed as a smooth submanifold) of $\mathbb R^n$. A diffeomorphism (i.e., a smooth invertible map) $\psi:\subs{\mathbb R^n} \rightarrow \subs{G}$, where $\subs{G}\subseteq G$, is called a \textit{local parametrization} of $G$. The map $\psi^{-1}$ is then a system of coordinates for $G$.
For instance, almost all (i.e., excepting a measure-zero subset) of $SO(3)$ can be parametrized using Euler angles.
If the support of $f$ is contained in $\subs{G}$, then the Euclidean mean of the random vector $\psi^{-1}(\tilde{\mf g})$ can be used to define a mean on $G$ \cite{zhou2021novel}. We refer to the quantity $\psi\left(\bmu_E(\psi^{-1}(\tilde{\mf g}))\right)$ as the \textit{parametric mean} of $\tilde{\mf g}$ w.r.t.~the parametrization $\psi$.

A special choice of parametrization is one that uses the $\exp$ map of the group. Let $\lbrace \mf E_i \rbrace_{i=1}^n$ be a basis for $\mathfrak g$.
Given a point $\mf h\in G$, define $\psi_{\mf h}(\mf x)\coloneqq \mf h \exp(x^i \mf E_i)$, where $\mf x = [\,x^1\ x^2\ \cdots\ x^n\,]\tran\in \mathbb R^n$ and the Einstein summation convention is used\footnote{That is, indices that appear twice within a term are summed, so that the expression $x^i \mf E_i$ is a notational shorthand for the sum $\sum_{i=1}^n x^i \mf E_i$ \cite[p. 18]{lee2012smooth}.}. Note that $\psi_{\mf h}(\mf 0_{n\times 1})=\mf h$;
the origin of $\mathbb R^n$ is mapped to $\mf h$. We let $(\,\cdot\,)^\vee:\mathfrak g \rightarrow \mathbb R^n$ denote the mapping from vectors in $\mathfrak g$ to their component-wise descriptions in the chosen basis. Given a vector $\mf Z\in\mathfrak g$, it is uniquely expressed as $\mf Z=z^i \mf E_i$ (where $z^i \in\mathbb R$), and we have $\mf Z^\vee = \begin{bmatrix}
    z^1 & \cdots & z^n
\end{bmatrix}\tran$.
From these definitions, it follows that $\psi_{\mf h}^{-1}(\mf g) = \log(\mf h^{-1}\mf g)^\vee$. We refer to $\psi_{\mf h}$ as the \textit{exponential parametrization} centered at $\mf h$.
The parametric mean defined using $\psi_{\mf h}$ is called the \textsl{log}-\textit{Euclidean} mean centered at $\mf h$, and has the expression\footnote{The evaluation of vector-valued integrals on $G$, such as the one in \cref{def:log-euclidean}, is explained in \Cref{app:mlg-calc}.}
\begin{align}\label{def:log-euclidean}
    \bmu_{\LE}(\tilde{\mf g};\mf h)
    &\coloneqq \mf h\exp\left(\int_{G'}  \log(\mf h^{-1}\mf g)\,
f(\mf g)\,d\mf  g \right) =\mf h\exp\left(\bmu_E(\log(\mf h^{-1}\tilde{\mf g}))\right).
\end{align}
For the exponential parametrization centered at $\mf I_m$, the identity element of $G$, the log-Euclidean mean becomes $\bmu_{\LE}(\tilde{\mf g};\,\mf I_m)=\exp\left(\bmu_E(\log(\tilde{\mf g}))\right)$ (cf. \cite[Ch. 6]{arsigny2006processing}).
Additionally, observe that when $\bmu\in \mu_G(\tilde{\mf g})$ is a group-theoretic mean, we have $\bmu_{\LE}(\tilde{\mf g};\bmu)=\bmu$.
Like the extrinsic Euclidean and projected means, parametric means can be evaluated in closed-form.
However, as we will show in \Cref{sec:properties,sec:examples}, parametric means generally lack many of the desirable properties that one expects of a mean of a random variable on $G$.

\section{Left, Right, and Bi-Invariant Distance Functions}
\label{sec:relationship}
The group-theoretic means $\mu_G(\tilde{\mf g})$ are defined without reference to an inner product or basis, while the Fr\'echet and Karcher means depend on the choice of distance function or Riemannian metric.
While there is a great degree of freedom when defining distances or metrics on $G$, one typically works with distances or metrics that are \textit{invariant} w.r.t.~the group operations. Invariant distances and metrics manifest as physically meaningful quantities in the applications of Lie groups to rigid-body dynamics \cite{vzefran1996choice,zefran1999metrics,chirikjian2015partial}, which motivates their use.
A distance function is said to be \textit{left-invariant} (resp., \textit{right-invariant}) if, for all $\mf g_1, \mf g_2, \mf h\in G$, 
$\mc D(\mf g_1, \mf g_2)$ $=$ $\mc D(\mf h \mf g_1,\mf h \mf g_2)$ (resp., $\mc D(\mf g_1, \mf g_2)=\mc D(\mf g_1\mf h, \mf g_2\mf h)$). For example, a \textit{bi-invariant} (i.e., both left and right-invariant) distance function for $SO(d)$ is $\mc D(\mf R, \mf Q)\coloneqq\lVert \mf R - \mf Q \rVert$, where $\mf R,\mf Q\in SO(d)$ and $d$ is a positive integer. Similarly, a left-invariant Riemannian metric is one whose geodesic distance function satisfies $\mc R(\mf g_1, \mf g_2)=\mc R(\mf h \mf g_1,\mf h \mf g_2)$.%
\footnote{
 See Appendix \ref{app:riemannian} for an equivalent characterization of left-invariant Riemannian metrics, which is stated in terms of the differential of the left multiplication map, $L_{\mf h}:\mf g \mapsto \mf h\mf g$.
}
While left- and right-invariant Riemannian metrics always exist on a given Lie group, there may not exist any bi-invariant Riemannian metrics; 
for example, $SE(d)$ does not admit any bi-invariant Riemannian metrics \cite[Prop. 7.2]{pennec2012exponential}.

To define an invariant Riemannian metric on $G$, one first defines an \textit{inner product} on the Lie algebra, $\mathfrak g$.
An inner product on $\mathfrak g$ is a symmetric bilinear map $\ip{\,\cdot\,}{\,\cdot\,}_{\mathfrak g}:\mathfrak g\times \mathfrak g \rightarrow \mathbb R$ satisfying the axioms in \cite[Sec. 5.1]{hornMatrixAnalysis2012}.
A choice of inner product on $\mathfrak g$ induces a \textit{norm} on $\mathfrak g$, given by $\lVert \mf X \rVert_{\mathfrak g} \coloneqq \sqrt{\ip{\mf X}{\mf X}}$. The inner product is said to be \textit{$\Ad$-invariant} if 
\begin{align}
\ip{\Ad_{\mf g} \mf X}{\Ad_{\mf g} \mf Y}&=\ip{\mf X}{\mf Y}
\end{align}
for all $\mf g\in G$ and $\mf X,\mf Y \in \mathfrak g$. 
As an example, consider the \textit{Frobenius inner product} for $\mathfrak{so}(d)$ defined as $\ip{\mf X}{\mf Y}_{\mathfrak{so}(d)}\coloneqq \tr(\mf X\tran \mf Y)$, which is $\Ad$-invariant.
However, if we consider a different inner product for $\mathfrak{so}(d)$, defined as $\ip{\mf X}{\mf Y}_{\mathfrak{so}(d)}^\prime \coloneqq \tr(\mf X\tran \mf Y \mf S)$ (where $\mf S$ is an arbitrary SPD matrix), such an inner product is \textit{not} $\Ad$-invariant in general. Similarly, if we extend the Frobenius inner product to all of $\mathfrak{gl}(d,\mathbb R)$ as $\ip{\mf X}{\mf Y}_{\mathfrak{gl}(d,\mathbb R)} \coloneqq \tr(\mf X\tran \mf Y)$, then this is {not} an $\Ad$-invariant inner product, since 
$$
\ip{\Ad_{\mf g}\mf X}{\Ad_{\mf g}\mf Y}_{\mathfrak{gl}(d,\mathbb R)} =
\tr\left(({\mf g}\mf X\mf g^{-1})\tran \mf g\mf Y\mf g^{-1}\right)\neq \tr(\mf X\tran \mf Y)$$ in general, where $\mf g\in GL(d,\mathbb R)$. 
We remark that the induced norm of the Frobenius inner product coincides with the Frobenius norm: $\lVert \mf X \rVert_{\mathfrak{gl}(d,\mathbb R)} = \lVert \mf X \rVert$.

Once an inner product is defined on $\mathfrak g$, it can be uniquely extended to define a {left-invariant} Riemannian metric on $G$, as detailed in \Cref{app:riemannian}. If the inner product on $\mathfrak g$ is $\Ad$-invariant, then the corresponding left-invariant Riemannian metric on $G$ will be bi-invariant.
The converse direction also holds, so that bi-invariant Riemannian metrics on $G$ and $\Ad$-invariant inner products on $\mathfrak g$ are in one-to-one correspondence.
Given a bi-invariant Riemannian metric on $G$ (equivalently, an $\Ad$-invariant inner product on $\mathfrak g$), the Riemannian log map is related to the Lie-theoretic $\log$ map as follows (cf. Lemma \ref{lem:bi-invariant-metrics}):
\begin{align}\rlog{\mf g}(\mf h) = \mf g \log(\mf g^{-1}\mf h)=  \log( \mf h\mf g^{-1})\mf g\quad \forall \mf g,\mf h\in G,
\label{eq:log-relationship}
\end{align}
where a capital `$\textrm{L}$' is used to distinguish the Riemannian log map from the Lie-theoretic one. Consequently, the geodesic distance function of a bi-invariant metric is $\mathcal R(\mf g, \mf h)=\lVert\log(\mf g^{-1}\mf h)\rVert_{\mathfrak g}$. For instance, the Frobenius inner product for $\mathfrak{so}(d)$ defines a bi-invariant Riemannian metric for $SO(d)$ whose geodesic distance function is $\mc R(\mf R, \mf Q)= \lVert\log(\mf R\tran\mf Q)\rVert$, where $\mf R,\mf Q\in SO(d)$.\footnote{In the literature, the inner product on $\mathfrak{so}(d)$ is sometimes defined as $\ip{\mf X}{\mf Y}_{\mathfrak{so}(d)}'\coloneqq\sfrac{1}{2}\,\tr(\mf X\tran \mf Y)$, such that the geodesic distance function incurs an additional factor of $\sfrac{1}{\sqrt 2}$ \cite{moakherMeansAveragingGroup2002}.}

In the absence of a bi-invariant Riemannian metric on $G$, the equality in (\ref{eq:log-relationship}) generally does not hold.\footnote{It is possible to endow the Lie group with its canonical \textit{Cartan-Schouten} connection (also called the $0$-connection), which defines a notion of a `geodesic' that satisfies a property analogous to \cref{eq:log-relationship}. This is the approach taken in \cite{pennec2013bi,pennec2020beyond}. However, these `geodesics' do not come with a corresponding notion of \textit{geodesic distance} unless the aforementioned connection is the Levi-Civita connection of some underlying Riemannian metric. Thus, the means defined in \cite{pennec2013bi,pennec2020beyond}, although stated using differential geometric concepts, are the same as the group-theoretic means considered in this article.}
For Lie groups like $SE(d)$, there do not exist \textit{any} Riemannian metrics for which (\ref{eq:log-relationship}) holds \cite{vzefran1996choice}, \cite[Lem. 3.1]{diatta2024cartan}. Even if $G$ does admit a bi-invariant metric, the application at hand may demand the use of a Riemannian metric for $G$ that is not bi-invariant.\footnote{For instance, $\tilde{\mf g} \in SO(3)$ can be used to describe the (uncertain) orientation of a rigid-body. In this case, a {weighted} inner product (i.e., one that differs from the Frobenius inner product) can be used to model the different moments of inertia of the rigid-body along its principal axes. The geodesics of the corresponding left-invariant Riemannian metric are not (in general) rotations about a fixed axis. Rather, they are the solutions to Euler's equations of rigid-body motion \cite[Sec. 1.2]{marsden2013introduction}.}
In what follows, we illustrate the relevance of bi-invariant Riemannian metrics to the problem of defining means on Lie groups.

 \begin{remark}[Bi-Invariant Riemannian Metrics]
A complete characterization of the connected Lie groups for which bi-invariant metrics exist is given in \cite[Lemma 7.5]{milnorCurvaturesLeftInvariant1976}: the Lie groups that admit bi-invariant metrics are those that are compact, Abelian (i.e., commutative), or a direct product of a compact Lie group and an Abelian Lie group. 
If $G$ is a simple Lie group, then the bi-invariant metric is unique up to scaling \cite[Lemma 7.6]{milnorCurvaturesLeftInvariant1976}. More generally, there are as many bi-invariant metrics as there are $\Ad$-invariant inner products on $\mathfrak g$ (cf. Lemma \ref{lem:bi-invariant-metrics}).
\end{remark}


Let $\lbrace \mf E_i \rbrace_{i=1}^n$ be a basis for $\mathfrak g$.
Any vector $\mf X\in \mathfrak g$ can be uniquely expressed as $\mf X=x^i\mf E_i$, where $x^i\in \mathbb R$ and the summation convention is used. The inner product on $\mathfrak g$ can be described using a SPD matrix $\mf W$, defined such that $\mf W_{ij}=\langle \mf E_i, \mf E_j \rangle$.\footnote{Alternatively, one way of \textit{defining} an inner product on $\mathfrak g$ is by choosing the basis $\lbrace \mf E_i \rbrace_{i=1}^n$ first, and then specifying the SPD matrix, $\mf W$.} Recalling the definition of the map $(\,\cdot\,)^\vee:\mathfrak g \rightarrow \mathbb R^n$ defined in \cref{sec:log-euclidean}, 
we observe that $$\lVert \mf X \rVert_{\mathfrak g} = \sqrt{\ip{\mf X}{\mf X}}_{\mathfrak g}=\sqrt{\ip{x^i\mf E_i}{x^j\mf E_j}_{\mathfrak g}}=\sqrt{\, {x^i}x^j \mf W_{ij}}=\lVert \mf X^\vee \rVert_{\mf W}.$$
The function $\mathcal L(\mf g, \mf h)\coloneqq \lVert\log(\mf g^{-1}\mf h)\rVert_{\mathfrak g}=\lVert \log(\mf g^{-1} \mf h)^\vee\rVert_{\mf W}$ satisfies some of the axioms of a distance function. In particular, we have $\mathcal L(\mf g, \mf h)=\mathcal L(\mf h, \mf g)$ and $\mathcal L(\mf g, \mf h)=0$ if and only if $\mf g=\mf h$.
Given these observations, we ask whether the group-theoretic means are critical points (i.e., local minimizers) of the following minimization problem:
\begin{align}
\underset{\mf h \in G}{\textrm{minimize}}\left(\int_{G'} \mathcal L(\mf g, \mf h)^2 f(\mf g) d \mf g\right).
\label{def:lie-karcher-mean}
\end{align}
That is, we ask whether a group-theoretic mean, like a Karcher mean, locally minimizes a cost function similar to \cref{eq:distance-mean-functional}. If the inner product is $\Ad$-invariant, then $\mathcal L$ is nothing but the geodesic distance function $\mc R$, so that (\ref{def:lie-karcher-mean}) indeed defines a Riemannian Fr\'echet mean whose local minimizers are $\mu_G(\tilde{\mf g})$ (i.e., $\mu_G(\tilde{\mf g})=\mu_K(\tilde{\mf g}, \mc R)$).
It turns out that $\Ad$-invariance of the inner product is not only sufficient, but also necessary for $\mu_G(\tilde{\mf g})$ to be critical points of (\ref{def:lie-karcher-mean}).

\vspace{5pt}
\begin{theorem}\label{thm:lie-karcher}
Let $\ip{\,\cdot\,}{\,\cdot\,}_{\mathfrak g}$ be an inner product on $\mathfrak g$.
\begin{itemize}
    \item If the inner product is $\Ad$-invariant, then given any random variable $\tilde{\mf g}\in G$, $\mu_G(\tilde {\mf g})$ is the set of critical points of (\ref{def:lie-karcher-mean}). 
    \item If the inner product is not $\Ad$-invariant, then there exists a random variable $\tilde{\mf g}\in G$ such that $\mu_G(\tilde{\mf g})$ are not critical points of (\ref{def:lie-karcher-mean}).
    \end{itemize}
\end{theorem}
\vspace{2pt}
\begin{proof}
Let $\mc C(\mf h)\coloneqq \int_{G'} \mc L(\mf g,\mf h)^2 f(\mf g) d \mf g$ be the cost function in \cref{def:lie-karcher-mean}.
By definition, $\bmu\in G$ is said to be a critical point of (\ref{def:lie-karcher-mean}) if we have $\ms X\mc C(\bmu)=0$
    for all $\mf X\in\mathfrak g$, where $\ms X \mc C$ denotes the Lie derivative of $\mc C$ along the left-invariant vector field generated by $\mf X$.
    As explained in \cref{app:mlg-calc}, $\ms X$ is a directional derivative operator whose direction is determined by $\mf X$.

    Consider the function $\psi_{\mf g}^{-1}(\mf h) \coloneqq \log(\mf g^{-1} \mf h)^\vee$ (which was introduced in \Cref{sec:log-euclidean}).
    From the definition of $\ms X$ and \cite[Sec. 5.5]{hall2013lie}, we have
    \begin{align}\label{eq:jacobian-log}
       \left(\ms X\psi_{\mf g}^{-1}\right)(\bmu) 
       &= \frac{d}{dt}\log(\mf g^{-1} \bmu \exp(t\mf X))^\vee\Big|_{t=0} = \mf J_{\log}(\mf g^{-1}\bmu)\,\mf X^\vee,
    \end{align}
    where, for $\mf h\in G$, the matrix $\mf J_{\log}(\mf h)$ represents the Jacobian of the $\log$ map at $\mf h$, and is given by the following series expansion (wherein the summation convention is not used):
    \begin{align}
    \mf J_{\log}(\mf h) = \mf I_n + \frac{1}{2}\bm\ad_{\log(\mf h)} + \sum_{i=1}^\infty \frac{\beta_{2i}}{(2i)!}\bm\ad_{\log(\mf h)}^{2i}.
    \end{align}
    Here, $\beta_i$ denotes the $i^{th}$ Bernoulli number and $\bm\ad_{\mf X}$ is the matrix of $\ad_{\mf X}$ in the basis $\lbrace \mf E_i \rbrace_{i=1}^n$, (cf. \Cref{app:mlg}). Using the chain rule and \cref{eq:jacobian-log}, we have
    \begin{align}
    \label{eq:lie-karcher-derivative-2}
       \ms X\mc C(\bmu) 
       &=\int_{G'} 2\,{\psi_{\mf g}^{-1}(\bmu)}\tran \mf W\,
       \left(\ms X\psi_{\mf g}^{-1}\right)(\bmu)\, f(\mf g)\, d \mf g \nonumber\\
       &= \int_{G'} 2\,{\log(\mf g^{-1} \bmu)^\vee}\tran \mf W\,\mf J_{\log}(\mf g^{-1}\bmu)\,\mf X^\vee f(\mf g)\, d \mf g. 
    \end{align} 
Let $\mf Z\coloneqq \log(\mf g^{-1}\bmu)$.
When the inner product is $\Ad$-invariant, it holds that $\bm\ad_{\mf Z}\tran \mf W = - \mf W\bm\ad_{\mf Z}$ (Lemma \ref{lem:bi-invariant-metrics}), so that 
$ \bm \ad_{\mf Z}^\top \mf W \mf Z^\vee=-\mf W [\mf Z,\mf Z]^\vee = \mf 0_{n\times 1}$. This fact (along with the fact that $\bmu\in\mu_G(\tilde{\mf g})$) makes the right-hand side of (\ref{eq:lie-karcher-derivative-2}) vanish for all $\mf X\in \mathfrak g$, showing one direction of the claim (i.e., \textit{sufficiency}). What remains to be shown is that $\Ad$-invariance is also \textit{necessary} for the right-hand side of (\ref{eq:lie-karcher-derivative-2}) to vanish in general.
    
Consider the random variable $\tilde{\mf g}$ that has the pdf $f(\mf g)=\delta(\mf h^{-1}\mf g)+\delta\big((\mf h^{-1})^{-1}\mf g\big),$
which is the empirical pdf consisting of two samples at $\mf h$ and $\mf h^{-1}$, where $\mf h\in G$. Since this pdf has the property that $f(\mf g)=f(\mf g^{-1})$ (i.e., it is symmetric about the identity element $\mf I_m\in G$), it readily follows that $\mf I_m \in \mu_G(\tilde{\mf g})$. On the other hand, we will show that $\mf I_m$ is (in general) not a critical point of (\ref{def:lie-karcher-mean}). Making the substitution $\bmu\mapsto \mf I_m$ in (\ref{eq:lie-karcher-derivative-2}), we get 
\begin{align*}
\ms X\mc C(\mf I_m)= 
 {\log(\mf h^{-1})^\vee}\tran\mf W \mf J_{\log}(\mf h^{-1})\mf X^\vee + {\log(\mf h)^\vee}\tran\mf W \mf J_{\log}(\mf h)\mf X^\vee= 0.\nonumber
\end{align*}
Since $\log\big(\mf h^{-1}\big)=-\log(\mf h)$, the terms in which $\log(\mf h)$ appears an odd number of times can be canceled to yield
\begin{align*}\label{eq:lie-karcher-derivative-3}
{\log(\mf h^{-1})^\vee}\tran\mf W\frac{1}{2}\bm\ad_{\log(\mf h^{-1})}\mf X^\vee + {\log(\mf h)^\vee}\tran\mf W\frac{1}{2}\bm\ad_{\log(\mf h)}\mf X^\vee={\log(\mf h)^\vee}\tran\mf W\bm\ad_{\log(\mf h)}\mf X^\vee=0.
\end{align*}
Since the preceding equation should hold true for all $\mf X\in\mathfrak g$, the condition for the point $\mf I_m$ to minimize (\ref{def:lie-karcher-mean}) is $\bm\ad_{\log(\mf h)}\tran \mf W\log(\mf h)^\vee=\mf 0_{n\times 1}$. Lemma \ref{lem:bi-invariant-metrics} shows that this holds for an arbitrary element $\mf h\in G$ if and only if the inner product is $\Ad$-invariant.
\end{proof}\vspace{4pt}

In addition to offering an alternative interpretation to the group-theoretic means and clarifying their relationship to the distance-based (i.e., Karcher and Fréchet) means, Theorem \ref{thm:lie-karcher} also offers a strategy to choose between the group-theoretic means in $\mu_G$. If $\mu_G(\mf g)$ has more than one element, then one can choose the element(s) of $\mu_G(\mf g)$ that minimize (\ref{def:lie-karcher-mean}) and discard the remaining elements (i.e., those that locally, but not globally minimize (\ref{def:lie-karcher-mean})). For example, it is shown in \cite[Fig. 4.1]{moakherMeansAveragingGroup2002} that, while a random variable on $SO(2)$ can have multiple group-theoretic means, it is often the case that there is a \textit{unique} group-theoretic mean 
that minimizes (\ref{def:lie-karcher-mean}).



\section{Properties of Means}
\label{sec:properties}
Each of the definitions of a \textit{mean} that were introduced in \Cref{sec:means} has distinctive properties that determine its suitability for a given statistical problem or application.
In this section, we discuss some of these properties, as well as the main considerations that are taken into account when choosing the most appropriate definition for a given application.

\subsection{Compatibility with Group Operations}\label{sec:invariance}
Given a (fixed) group element $\mf h\in G$, we can define the left-shifted random variable, $\mf h \tilde {\mf g}$. The pdf of $\mf h \tilde {\mf g}$ is the function $\mf g \mapsto f(\mf h^{-1} \mf g)$.\footnote{Here, it is assumed that the pdfs of $\tilde {\mf g}$ and $\mf h \tilde {\mf g}$ are each defined w.r.t.~the left Haar measure (c.f. \Cref{app:mlg-calc}).}
If $\bmu$ is a mean of $\tilde {\mf g}$, one can ask whether $\mf h\bmu$ is a mean of $\mf h \tilde {\mf g}$. This property is called \textit{left-invariance}, and it holds for group-theoretic means. To see this, observe that
\begin{align*}
\int_{G'} \log(\bmu^{-1} \mf g) f(\mf g) d \mf g &= \int_{G'} \log(\bmu^{-1} \mf h^{-1} \mf g) f(\mf h^{-1} \mf g) d \mf g = \int_{G'} \log((\mf h\bmu)^{-1} \mf g) f(\mf h^{-1} \mf g) d \mf g.
\end{align*}
Hence, if $\bmu\in\mu_G(\tilde {\mf g})$, then $\mf h\bmu \in \mu_G(\mf h \tilde {\mf g})$. Succinctly, we can write $\mu_G(\mf h \tilde {\mf g})=\mf h\bmu_G(\tilde {\mf g})$ for all $\mf h\in G$. 
When $G=(\mathbb R^n, +)$ and $\tilde {\mf x}$ is an $\mathbb R^n$-valued random variable, left-invariance of $\mu_G$ reduces to the well-known property: $\mu_G(\tilde{\mf x}-\mf h)=\mu_G(\tilde{\mf x})-\mf h$, where $\mf h\in \mathbb R^n$.

The Lie group $G$ is said to be  \textit{unimodular} if the left Haar measure is (in addition to being left-invariant) right-invariant. Unimodularity (i.e., the existence of a bi-invariant Haar measure) is a weaker condition than the existence of a bi-invariant metric. In particular, $SE(d)$ and $GL(d, \mathbb R)$ are unimodular, which can be verified using the conditions for unimodularity that are presented in \Cref{app:mlg-calc}. When $G$ is unimodular, the pdf of $\tilde {\mf g} \mf h$ is $f(\mf g \mf h^{-1})$. From (\ref{def:group-theoretic-mean-right}) and the invariance properties of the Haar measure, it follows that the group-theoretic means are right-invariant as well.
For this reason, \cite{pennec2012exponential,pennec2020beyond} refers to the group-theoretic means as \textit{bi-invariant means}.

The extrinsic Euclidean mean satisfies
$
\bmu_E(\mf A\tilde{\mf g})=\mf A\bmu_E(\tilde{\mf g})$ {and} $\bmu_E(\tilde{\mf g}\mf A)=\bmu_E(\tilde{\mf g})\mf A$
for any matrix $\mf A\in \mathbb R^{m\times m}$, which implies that it is bi-invariant.
The projected mean is left- (resp., right-) equivariant if $\mathcal P$ commutes with the left- (resp., right-) multiplication map of $G$. For example, if $\mathcal P(\mf h \mf A)= \mf h \mathcal P(\mf A)$ for all $\mf A\in \mathbb R^{m\times m}$ and $\mf h\in G$, then we can show that the projected mean is left-invariant:\begin{align}
\mathcal P(\bmu_E(\mf h \tilde {\mf g}))=\mathcal P(\mf h\bmu_E(\tilde {\mf g}))=\mf h \mathcal P(\bmu_E(\tilde {\mf g})).
\end{align}
For example, it is readily verified that the map $\mc P_{SO}$ introduced in \Cref{sec:extrinsic}
is bi-invariant, and therefore defines a bi-invariant projected mean for $SO(d)$.
A set of Fréchet means (defined with respect to a particular choice of distance function) is left-, right-, or bi-invariant if the corresponding distance function has the desired invariance property.  
Similarly, Karcher means are left-, right-, or bi-invariant if the Riemannian metric used to define them has the corresponding invariance \cite{pennec2020beyond}, but the same does not hold for a general Riemannian metric. 


\subsection{Defining a Notion of Covariance}
One of the main applications of the concept of a mean is to summarize the random variable $\tilde {\mf g}$. In this context, it is called a \textit{first-order statistic}. While a first-order statistic describes the `center of mass' of a pdf, a \textit{second-order statistic} characterizes its dispersion/spread about the center of mass. In the case of the Euclidean mean, the most commonly used second-order statistics are the variance and the covariance.

The $\vecn(\,\cdot\,)$ operation introduced in \Cref{sec:euclidean}, can be used to define an $m^2 \times m^2$ matrix which we refer to as the \textit{extrinsic Euclidean covariance} matrix:
\begin{align}
    \Cov_E(\tilde {\mf g}) \coloneqq \int_{G}{\vecn\big(\,\mf g - \bmu_E(\tilde{\mf g})\big)\,\vecn\big(\,\mf g - \bmu_E(\tilde{\mf g})\big)\tran}\,f(\mf g)\,d \mf g.
    \label{eq:euclidean-covariance}
\end{align}
An extrinsic covariance for the projected mean can be defined similarly, by replacing $\bmu_E(\tilde{\mf g})$ with $\mc P(\bmu_E(\tilde{\mf g}))$. In either case, the dimension of the covariance matrix may far exceed the intrinsic dimension of the Lie group.
The trace of the covariance can be used to define the \textit{Euclidean variance}:
\begin{align}
    \textrm{var}_E(\tilde {\mf g}) &\coloneqq {\tr\big(\Cov_E(\tilde {\mf g})\big)}
    ={\int_{G}\left\lVert\vecn\big(\tilde {\mf g}-\bmu_E(\tilde{\mf g})\big)\right\rVert^2 f(\mf g)\,d \mf g}={\int_{G}\lVert\tilde {\mf g}-\bmu_E(\tilde{\mf g})\rVert^2 f(\mf g)\,d \mf g}.
    \label{eq:euclidean-variance}
\end{align}

If $\bmu\in\mu_F(\tilde {\mf g};\mathcal D)$ is a Fréchet mean, one can define the \textit{Fréchet variance} of $\tilde {\mf g}$ as the minimal value of the cost function (\ref{eq:distance-mean-functional}), which is attained at $\bmu$:
\begin{align}
\textrm{var}_F(\tilde {\mf g};\mathcal D) &\coloneqq 
    {\int_{G'}\mathcal D(\bmu, \mf g)^2 f(\mf g)\,d \mf g},
\end{align}
where the choice of $\bmu$ in $\mu_F(\tilde {\mf g};\mathcal D)$ does not matter.
However, since the definition of a Fréchet mean does not (in general) involve vector operations, there is no natural notion of a \textit{covariance}; one can describe the spread of $\tilde {\mf g}$ about its mean using $\textrm{var}_F$, but the directional aspects of this spread cannot be captured using $\textrm{var}_F$.

Given the basis $\lbrace \mf E_i \rbrace_{i=1}^n$, the group-theoretic covariance and variance are defined as \cite{barfoot2014associating,wang2008nonparametric}
\begin{align}
\label{def:covariance-group}
\bm\Sigma_G(\tilde {\mf g};\bmu) &\coloneqq \int_{G'} \log(\bmu^{-1}\mf g)^\vee{\log(\bmu^{-1}\mf g)^\vee}\tran f(\mf g) \,d \mf g \\
\label{def:variance-group}
\textrm{var}_G(\tilde {\mf g};\bmu) &\coloneqq {\tr\big(\bm\Sigma_G(\tilde {\mf g};\bmu)\big)} = {\int_{G'}\big\lVert \log(\bmu^{-1} \mf g)^\vee\big\rVert^2\,f(\mf g) \,d \mf g},
\end{align}
where $\bmu\in \mu_G(\tilde {\mf g})$.
While the definition of $\mu_G$ did not depend on the choice of basis for $\mathfrak g$, the definitions in \cref{def:covariance-group} and \cref{def:variance-group} do depend on the choice of basis. In particular, if we choose an orthonormal basis for $\mathfrak g$ (such that $\mf W=\mf I_n$), then we have
\begin{align}
\textrm{var}_G(\tilde {\mf g};\bmu)
= {\int_{G'}\lVert \log(\bmu^{-1} \mf g)\rVert_{\mathfrak g}^2\,f(\mf g) \,d \mf g},\end{align} 
which is precisely the cost function that was considered in \Cref{sec:relationship}. Thus, if the inner product on $\mathfrak g$ is $\Ad$-invariant, and if one defines the covariance of $\tilde {\mf g}$ using an orthonormal basis, then the group-theoretic means are minimizers of the map $\bmu \mapsto \textrm{var}_G(\tilde {\mf g};\bmu)$.

The Karcher covariance and variance can be defined similarly, by letting $\bmu \in \mu_K(\tilde {\mf g}; \mathcal R)$ and using $\bmu^{-1}\rlog{\bmu}(\mf g)$ in place of $\log(\bmu^{-1} \mf g)$ in (\ref{def:covariance-group}). 
If the inner product on $\mathfrak g$ is $\Ad$-invariant and $\lbrace \mf E_i \rbrace_{i=1}^n$ is orthonormal, then the Karcher covariance (defined w.r.t.~the corresponding bi-invariant metric) is the same as the group-theoretic covariance. 
When using a parametric mean (as defined in \Cref{sec:log-euclidean}), the Euclidean covariance of the vector-valued random variable, $\psi^{-1}(\tilde {\mf g})$, can be used as a measure of dispersion.


\subsection{Compatibility with Symmetries}\label{sec:symmetric}
Given a point $\bmu \in G$, we say that the pdf $f$ is \textit{symmetric about $\bmu$} if it has the property $f(\bmu\mf g)=f(\bmu \mf g^{-1})$ (equivalently, $f(\mf g \bmu)=f(\mf g^{-1}\bmu)$) for all $\mf g\in G$. Let $\rho$ denote the pdf of the shifted random variable $\bmu^{-1}\tilde{\mf g}$, such that $\rho(\mf g)=f(\bmu\mf g)$. It follows that $f$ is symmetric about $\bmu$ if and only if $\rho$ is symmetric about the identity element, i.e., $\rho(\mf g)=\rho(\mf g^{-1})$.
As the forthcoming theorem shows, if $G$ is unimodular, then the symmetry of $f$ about $\bmu$ implies that $\bmu\in \mu_G(\tilde {\mf g})$. Furthermore, the group-theoretic means minimize a cost function similar to (\ref{def:lie-karcher-mean}) wherein the norm on $\mathfrak g$ is defined in terms of the group-theoretic covariance. Unlike the analysis in \Cref{sec:relationship}, the analysis here does not assume that the inner product on $\mathfrak g$ is $\Ad$-invariant.
To state the result, we need the following lemma.

\begin{lemma}\label{lem:anti-symmetric}
Let $\bm\upeta$ be a vector-valued function on $G$ that satisfies $\bm\upeta (\mf h \mf g)= -\bm\upeta (\mf h \mf g^{-1})$ for all $\mf g \in G$, where $\mf h\in G$ is some (fixed) point. If $G$ is unimodular, then
$
\int_G\bm\upeta(\mf g)d\mf g = \mf 0
$.
\end{lemma}
\begin{proof}
Since $G$ is unimodular, the Haar integral is invariant under left multiplication as well as the inversion map $\mf g \mapsto \mf g^{-1}$. Using these facts, we see that 
\begin{align*}
\int_G\bm\upeta(\mf g)d\mf g = \int_G\bm\upeta(\mf h\mf g)d\mf g = -\int_G\bm\upeta(\mf h\mf g^{-1})d\mf g= -\int_G\bm\upeta(\mf h\mf g) d\mf g= - \int_G\bm\upeta(\mf g)d\mf g = \mf 0.
\end{align*}
\end{proof}

\vspace{3pt}
\begin{theorem}\label{thm:symmetric}
    If $G$ is unimodular and $f$ is symmetric about $\bmu$, then $\bmu\in \mu_G(\tilde {\mf g})$. Furthermore, if $\bm\Sigma_G(\tilde {\mf g}; \bmu)$ is invertible, then $\bmu$ is a critical point of the problem
    \begin{align}\label{def:lie-karcher-mean-weighted}
    \underset{\mf h \in G}{\textrm{minimize}}\left(\int_{G'} \lVert \log(\mf h^{-1} \mf g)^\vee \rVert^2_{\bm\Sigma^{-1}} f(\mf g) d \mf g\right),
    \end{align}
    where $\bm\Sigma\coloneqq {\bm\Sigma_G(\tilde {\mf g}; \bmu)}$.
\end{theorem}
\begin{proof}
First, we observe that the function $\bm\upeta(\mf g)\coloneqq \log(\bmu^{-1} \mf g)$ satisfies $\bm\upeta (\bmu \mf g)= -\bm\upeta(\bmu\mf g^{-1})$ (since $\log(\mf g)=-\log(\mf g^{-1})$),
whereupon \cref{lem:anti-symmetric} is used to conclude that $\bmu\in \mu_G(\tilde {\mf g})$.
To show the second part of the theorem, we proceed similarly to the proof of \cref{thm:lie-karcher}. The condition for minimization of (\ref{def:lie-karcher-mean-weighted}) is
    \begin{align}\nonumber
\int_{G'} \,{\bm\upeta(\mf g)^\vee}\tran \bm\Sigma^{-1} \mf J_{\log}({\bm\upeta(\mf g)})\mf X^\vee\, f(\mf g)\, d \mf g=\frac{1}{2}\int_{G'} \,{\bm\upeta(\mf g)^\vee}\tran \bm\Sigma^{-1}\bm\ad_{\bm\upeta(\mf g)}\mf X^\vee\, f(\mf g)\, d \mf g = 0
\end{align}
$\forall \mf X\in \mathfrak g$, where we used the fact that each of the terms of $\mf J_{\log}({\bm\upeta(\mf g)})$ in which $\bm\upeta(\mf g)$ appears an odd number of times vanish due to \cref{lem:anti-symmetric}. 
Expressing $\bm\upeta(\mf g)$ as $\eta^i(\mf g)\mf E_i$,
the condition for minimization can be written using components, as
\begin{align}
\int_{G'} \,\eta^i(\mf g) \Sigma^{ij}
\eta^k(\mf g)  C_{k\ell}^j\, f(\mf g)\, d \mf g =0\quad\forall \ell\in\lbrace 1,\ldots,n\rbrace,
\label{eq:sym-proof-cond1}
\end{align}
where $\Sigma^{ij}$ is the $(i,j)^{th}$ element of $\bm\Sigma^{-1}$ and $\lbrace C_{ij}^k\rbrace_{i,j,k=1}^n$ are the \textit{structure constants} of $\mathfrak g$; see \Cref{app:mlg} for details. Let $\Sigma_{ij}=\int_{G'}\eta^i(\mf g)\eta^j(\mf g) f(\mf g) d \mf g$ be the $(i,j)^{th}$ component of $\bm\Sigma$. We can rewrite \cref{eq:sym-proof-cond1} as
\begin{align}
\Sigma^{ij}C_{k\ell}^j\int_{G'} \,\eta^i(\mf g) 
\eta^k(\mf g)  \, f(\mf g)\, d \mf g =\Sigma^{ij}C_{k\ell}^j\Sigma_{ki} =C_{k\ell}^j \delta_{jk} =C_{k\ell}^k = 0
\label{eq:symmetric-unimodularity}
\end{align}
$\forall \ell\in\lbrace 1,\ldots,n\rbrace$. Here, $\delta_{ij}$ is the $(i,j)^{th}$ element of $\mf I_n$.
The last equality in (\ref{eq:symmetric-unimodularity}) is precisely the condition for the unimodularity of $G$ (c.f. App. \ref{app:mlg-calc}). Hence, the conditions for $\bmu$ to be a critical point of (\ref{def:lie-karcher-mean-weighted}) are that $G$ is unimodular and $\bm\Sigma$ is invertible.
\end{proof}\vspace{2pt}
The following corollary follows readily from the proof of \cref{thm:symmetric}.

\vspace{3pt}
\begin{corollary}\label{cor:mean_identity}
If $\bmu\in\mu_G(\tilde {\mf g})$ and if the pdf $f$ is sufficiently concentrated (so that terms on the order of $O(\lVert\log(\bmu^{-1} \mf g)^\vee\rVert^3)$ may be neglected), then $\bmu$ is a critical point of (\ref{def:lie-karcher-mean-weighted}). 
\end{corollary}

\noindent
Therefore, a group-theoretic mean $\bmu\in \mu_G(\tilde {\mf g})$ will locally minimize the cost function given in \cref{def:lie-karcher-mean} (which is defined w.r.t.~an inner product on $\mathfrak g$) if either of the following conditions hold:
\begin{enumerate}[(i)]
\item (\cref{thm:lie-karcher}) the inner product on $\mathfrak g$ is $\Ad$-invariant, or 
\item (\cref{thm:symmetric}) $\bm\Sigma$ is invertible, $f$ is symmetric about $\bmu$, and the inner product is defined such that $\mf W=\bm\Sigma^{-1}$,
\end{enumerate}
where $\bm\Sigma\coloneqq\bm\Sigma_G(\tilde {\mf g};\bmu)$.


\subsection{Compatibility with Convolution}
Given a pair of independent random variables $\tilde{\mf g}_1$ and $\tilde{\mf g}_2$ in $G$, we consider the problem of characterizing their group-theoretic product, $\tilde{\mf g}_1\tilde{\mf g}_2$. 
Let the pdfs of $\tilde {\mf g}_1$ and $\tilde {\mf g}_2$ (as defined w.r.t. the left Haar measure) be $f_1$ and $f_2$, respectively. When $G$ is unimodular, the pdf of $\tilde {\mf g}_1 \tilde {\mf g}_2$ is given by the group-theoretic \textit{convolution} of $f_1$ and $f_2$, defined as follows \cite{chirikjian1998numerical,chirikjian2000engineering}:
\begin{align}\label{def:convolution}
\big(f_{1}\conv f_2\big)(\mf h)&\coloneqq\int_{G} f_1(\mf g)f_2(\mf g^{-1}\mf h)d \mf g.
\end{align}
Due to the computational difficulty of evaluating the integral in (\ref{def:convolution}), one is often interested in computing the mean and covariance of $\tilde {\mf g}_1 \tilde {\mf g}_2$ as a function of the means and covariances of $\tilde {\mf g}_1$ and $\tilde {\mf g}_2$. This bypasses the computationally difficult step of evaluating $f_1 \conv f_2$, and is amenable to real-time applications where the mean and covariance of $\tilde {\mf g}_1 \tilde {\mf g}_2$ may be needed for time-sensitive tasks.

In the case of $G=(\mathbb R^n, +)$, the group-theoretic product of random variables $\tilde{\mf x}_1$ and $\tilde{\mf x}_2$ is their vector sum $\tilde{\mf x}_1+\tilde{\mf x}_2$, whose mean and covariance can be described in terms of the mean and covariance of $\tilde{\mf x}_1$ and $\tilde{\mf x}_2$ as
\begin{align*}
\bmu_E(\tilde{\mf x}_1 + \tilde{\mf x}_2) = \bmu_E(\tilde{\mf x}_1) + \bmu_E(\tilde{\mf x}_2),\quad\text{and}\quad \bm\Sigma_{E}(\tilde{\mf x}_1 + \tilde{\mf x}_2) = \bm\Sigma_{E}(\tilde{\mf x}_1) + \bm\Sigma_{E}(\tilde{\mf x}_2),
\end{align*}
respectively.
The group-theoretic means of $\tilde {\mf g}_1 \tilde {\mf g}_2$ satisfy this property approximately, in the following sense. Let $\bmu_1\in \mu_G(\tilde {\mf g}_1)$ and $\bmu_2\in \mu_G(\tilde {\mf g}_2)$. If the pdf is sufficiently concentrated (such that third-order terms can be neglected), then it holds that $\bmu_1 \bmu_2 \in \mu_G(\tilde {\mf g}_1 \tilde {\mf g}_2)$. The group-theoretic covariance of $\tilde {\mf g}_1 \tilde {\mf g}_2$ about this mean is 
$$
\bm\Sigma_{12} \approx \bm \Ad_{\bmu_2}^{-1} \bm \Sigma_1 {\bm\Ad_{\bmu_2}^{-1}}\tran + \bm \Sigma_2,$$ 
where $\bm \Sigma_1\coloneqq \bm\Sigma_G(\tilde {\mf g}_1; \bmu_1)$, $\bm \Sigma_2\coloneqq \bm\Sigma_G(\tilde {\mf g}_2; \bmu_2)$, and $\bm\Sigma_{12} \coloneqq \bm\Sigma_G(\,\tilde {\mf g}_1 \tilde {\mf g}_2\,; \bmu_1\bmu_2)$.
Higher-order correction terms can be introduced when a greater degree of accuracy is required \cite{wang2008nonparametric,park2010path}.

\begin{remark}
    One can also consider a stochastic process on $G$, which can be used to describe the continuous-time evolution of uncertainty in a random variable. Suppose a stochastic process $(\tilde {\mf g}_t)_{t\geq 0}$ on $G$ is defined via a \textit{stochastic differential equation (SDE)}. The corresponding family of pdfs $(f_t)_{t\geq 0}$ (where $f_t$ is the pdf of $\tilde{\mf g}_t$) satisfies a partial differential equation, known as the \textit{Fokker-Planck equation} \cite{chirikjian2011stochastic}. While it is generally infeasible to solve the Fokker-Planck equation explicitly, it is possible to update the group-theoretic mean and covariance of $\tilde{\mf g}_t$ using ordinary differential equations \cite{park2010path,ye2024uncertainty}.
\end{remark}

Analogous properties hold for the Euclidean mean as well. Using the invariance properties of the Haar integral, we have
\begin{align}
\bmu_E(\tilde {\mf g}_1 \tilde {\mf g}_2) &= 
\int_G\int_G \mf g_1 \mf g_2 \,f_1(\mf g_1)f_2(\mf g_2)\,d \mf g_1 \,d\mf g_2 
= \bmu_E({\tilde {\mf g}_1})\,\bmu_E({\tilde {\mf g}_2}).
\end{align}
It is clear that the projected mean does not have an analogous property, since $\mc P$ is generally a nonlinear map that does not commute with integration. In general, the parametric, Fr\'echet, and Karcher means (defined w.r.t.~an arbitrary Riemannian metric) of $\tilde{\mf g}_1\tilde{\mf g}_2$ are also unrelated to the corresponding means of $\tilde{\mf g}_1$ and $\tilde{\mf g}_2$.

\subsection{Computational Feasability}\label{sec:computation}
For some Lie groups, such as $SO(d)$ and $SE(d)$ ($d=2$ or $3$), the group-theoretic $\exp$ and $\log$ maps have closed-form expressions \cite{chirikjian2011stochastic,barfoot2024state}. More generally, the matrix exponential and logarithm functions must be used to evaluate $\exp$ and $\log$, respectively; these functions are implemented in most linear algebra toolboxes. On the other hand, the Riemannian exponential map is generally computed by solving ordinary differential equations (ODEs). While the geodesics of an arbitrary Riemannian metric are solutions to a system of second-order ODEs known as the \textit{geodesic equation}, the geodesics of a left-invariant Riemannian metric are solutions of first-order ODEs (c.f. \Cref{app:geodesics}), which simplifies their implementation. In either case, computational tools such as \verb|geomstats| (for Python) \cite{guigui2023introduction} and \verb|Manifolds.jl| (for Julia) \cite{manifolds.jl} can be used to evaluate the Riemannian exponential and log maps on Lie groups. Nevertheless, if a bi-invariant metric for $G$ exists, then one prefers to work with the bi-invariant metric due to the computational ease of implementing the Lie-theoretic $\exp$ and $\log$ maps.

Consider the map $\bmu \mapsto \bmu_{\LE}(\tilde {\mf g};\bmu)$, where $\bmu_{\LE}$ is the log-Euclidean mean. As discussed in \Cref{sec:log-euclidean}, a fixed point of this map (i.e., a point $\bmu^\star$ that satisfies $\bmu^\star=\bmu_{\LE}(\tilde {\mf g};\bmu^\star)$) 
is a group-theoretic mean. Hence, the fixed point iteration $\bmu_{k+1}=\bmu_{\LE}(\tilde {\mf g};\bmu_k)$ can be used for computing the group-theoretic means, as noted in \cite{pennec2013bi,wang2008nonparametric}. 
A slight modification of this algorithm (that replaces the Lie-theoretic exponential and log maps with the Riemannian ones) yields a procedure for computing Karcher means, albeit with the added computational cost of solving the geodesic equation. 
Fr\'echet means can be computed using optimization algorithms on Lie groups, such as gradient descent.
There also exist second-order algorithms for computing group-theoretic and Fr\'echet means \cite{owren2000newton,argyros2009newton}.



\section{Examples of Lie Groups}
\label{sec:examples}
In this section, we consider some of the prototypical examples of finite-dimensional Lie groups that typically arise in engineering and scientific applications. 



\subsection{$SO(2)$ and $U(1)$}
\label{sec:SO2}
The \textit{Orthogonal} group $O(2)$ is the set of $2\times 2$ orthogonal matrices, such that $\mf R\in O(2) \Leftrightarrow \mf R^{-1}=\mf R\tran$. However, when viewed as a manifold, $O(2)$ has two disconnected components -- there exist pairs of elements $\mf R,\mf Q\in O(2)$ with $\det(\mf R)=+1$ and $\det(\mf Q)=-1$ that lie in different connected components of $O(2)$, such that there is no geodesic in $O(2)$ connecting $\mf R$ and $\mf Q$. Since $\det(\mf R^{-1}\mf Q)=-1$, the matrix $\mf R^{-1}\mf Q$ is outside the domain of the group-theoretic $\log$ map as well. Due to these difficulties, one typically restricts to the \textit{Special Orthogonal} group $SO(2)$, which comprises the elements of $O(2)$ whose determinant is $+1$, to ensure that the group-theoretic and Karcher means are well-defined.

$SO(2)$ is isomorphic to the \textit{Unitary} group $U(1)$, which can be viewed as the set of unitary matrices in ${GL}(1,\mathbb C)$, i.e., the unit complex numbers. That is, the anti-clockwise rotation by $\theta$ radians can be represented as the unit complex number $e^{\ii\theta} \in U(1)$, or as the $2\times 2$ rotation matrix
\begin{align}
\mf R(\theta) = \begin{bmatrix}
\cos\theta & -\sin\theta \\
\sin\theta & \cos\theta
\end{bmatrix} \in SO(2).
\label{eq:so2-param}
\end{align}
In this sense, $SO(2)$ and $U(1)$ can each be identified with the unit circle of $\mathbb R^2$, denoted as $S^1$ (i.e., the one-dimensional sphere).
For this reason, the theory of random variables on these Lie groups is often referred to as \textit{circular} statistics. 
Since points on the circle correspond bijectively to the \textit{Direction of Arrival (DoA)} of a wireless signal being received at a sensor array, the uncertainty in the DOA can be modeled using circular statistics \cite{jacob2013DOA,wozniak2019self}.
The identification of $SO(2)$ with the unit circle is also useful for defining pdfs on the circle, such as the wrapped Gaussian, wrapped Cauchy, von Mises, Fisher, and Bingham distributions \cite{jupp1989unified,jammalamadaka2001topics,mardia2009directional}. 

Given a random variable $\tilde{\mf R}$ in $SO(2)$, the matrix $\bmu_E(\tilde{\mf R})$ is generally not an element of $SO(2)$. Nevertheless, we can assume that it is invertible (if not, an arbitrarily small amount of noise can be added to the data to ensure its invertibility). Given an invertible matrix $\mf A$, we can consider the map $\mathcal P_{SO}(\mf A) \coloneqq (\mf A \mf A\tran)^{-\frac{1}{2}}\mf A=\mf A (\mf A\tran \mf A)^{-\frac{1}{2}}$ that extracts the rotation matrix in the polar decomposition of $\mf A$. This defines a bi-invariant projected mean for $SO(2)$ that is frequently used in the literature \cite{moakherMeansAveragingGroup2002,sarlette2009consensus}. Additional care needs to be taken to ensure that the projected Euclidean mean has a positive determinant, e.g., by multiplying $\mathcal P_{SO}(\mf A)$ with $\mathrm{sign}(\det(\mf A))$. 


To define the Fr\'echet mean on $SO(2)$, one can consider the {Euclidean distance} between points in $SO(2)$, which is the function $\mathcal D(\mf R, \mf Q)\coloneqq\lVert \mf R - \mf Q \rVert$. 
The resulting Fr\'echet mean is called the \textit{induced arithmetic mean} in \cite{sarlette2009consensus}, and is shown to be the same as the projected mean defined via $\mathcal P_{SO}$ above. The Riemannian Fr\'echet mean is defined using the geodesic distance between points in $SO(2)$, which (for the Frobenius inner product introduced in \Cref{sec:relationship})
is $\mc R(\mf R,\mf Q)=\lVert\log(\mf R\tran \mf Q)\rVert$.\footnote{In fact, since $SO(2)$ is Abelian, $\Ad_{\mf R}$ is the identity map on $\mathfrak g$ for all $\mf R\in SO(2)$. It follows that any inner product on $\mathfrak {so}(2)$ is $\Ad$-invariant. A choice of inner product on $\mathfrak {so}(2)$ amounts to a choice of units for measuring angles; since $\mathfrak {so}(2)$ is one-dimensional, the weighting ``matrix" $\mf W$ is in this case a positive number.} It follows that the group-theoretic means are critical points of the Riemannian Fr\'echet mean defined by $\mathcal R$. They are not unique; if the pdf of $\tilde{\mf R}$ is symmetric about the identity element (i.e., $f(\mf R)=f(\mf R\tran)$), it can be shown that $\mf I_2$ and $-\mf I_2$ are both in $\mu_G(\tilde{\mf R})$. Nevertheless, in general only one of these means will minimize the cost function (\ref{def:lie-karcher-mean}). As observed in \cite{moakherMeansAveragingGroup2002}, this offers a strategy to compute the Fr\'echet mean on $SO(2)$; one first computes $\mu_G(\tilde {\mf g})$ and then selects the mean(s) that minimize (\ref{def:lie-karcher-mean}), which may be unique.

\begin{remark}
If we identify $SO(2)$ with $U(1)$ via the map $\mf R(\theta)\mapsto e^{\ii \theta}$, then $\mc D(\mf R, \mf Q)$ is nothing but the \textit{chordal distance} (i.e., the length of the chord) between $\mf R$ and $\mf Q$, whereas $\mc R(\mf R, \mf Q)$ is the arclength between $\mf R$ and $\mf Q$.
\end{remark}


\subsection{$SO(3)$ and $SU(2)$}
\label{sec:SU2}
The \textit{Special Unitary} Lie group $SU(2)$ is isomorphic to the group of unit quaternions. The unit quaternions (and therefore, $SU(2)$) can be visualized as the unit sphere $S^3$ of $\mathbb R^4$.
Each unit quaternion can be mapped to a matrix in $SO(3)$, due to which quaternions are used extensively in fields like computer graphics and spaceflight dynamics \cite{psiaki2006estimation,markleyAveragingQuaternions2007,vince2011quaternions} to encode orientations of rigid bodies.
Since the quaternions $\mf q$ and $-\mf q$ both map to the same rotation, we say that $SU(2)$ is a \textit{double cover} of $SO(3)$. While $SO(3)$ and $SU(2)$ have isomorphic Lie algebras (i.e., $\mathfrak{so}(3)\cong \mathfrak{su}(2)$), they are not isomorphic as Lie groups. 

Due to the identification of $SU(2)$ with $S^3$, the existing literature on \textit{spherical statistics} \cite{buss2001spherical}, which deals with pdfs on spheres, can be applied to the study of statistical distributions on $SU(2)$. In particular, the function $\mc P_{S^3}: \mf x \mapsto \mf x/\lVert \mf x \rVert$ can be used to project points from $\mathbb R^4\backslash \lbrace 0 \rbrace$ onto $S^3$. Similar to circular statistics, in spherical statistics, the Fisher and Bingham distributions can be defined as pdfs on $S^3$.
However, when describing $SO(3)$-valued random variables, one should ensure that the pdf thus defined on $S^3$ is antipodally symmetric (as in the case of the Bingham distribution). That is, the pdf on $S^3$ should have the same value at $\mf q$ and $-\mf q$, as these points (viewed as quaternions) correspond to the same element in $SO(3)$.

The projection function also provides a way to define the mean of quaternions. However, the projected mean of quaternions has undesirable features when the quaternions are meant to represent spatial rotations (i.e., elements of $SO(3)$). While points on the quaternion sphere that are (almost exactly) antipodal correspond to the same rotation, their Euclidean mean is close to the origin of $\mathbb R^4$ the projected mean is highly sensitive to perturbations of either point \cite{markleyAveragingQuaternions2007}.
Nevertheless, the quaternion-based method of averaging rotations remains popular in the state estimation literature, which is typically focussed on concentrated distributions (i.e., pdfs of random variables that have a small variance) \cite{psiaki2006estimation}.

In \cite{markleyAveragingQuaternions2007},
the authors propose to use a Fr\'echet mean on $SO(3)$ rather than the projected mean of quaternions. Specifically, let $\mu_{F}(\tilde {\mf g}; \mathcal D)$ be the Fr\'echet mean defined by the Euclidean distance function $\mathcal D(\mf R, \mf Q) \coloneqq \lVert \mf R - \mf Q \rVert$. The Fr\'echet mean cost function \cref{eq:distance-mean-functional} is then the same as the cost function in \cref{eq:euclidean-matrix-mean-frob}, which is minimized by the Euclidean mean $\mu_E(\tilde {\mf R})$. However, unlike $\mu_E(\tilde {\mf R})$, the Fr\'echet means in $\mu_F(\tilde {\mf g}; \mathcal D)$ are constrained to lie in $SO(3)$. It is shown in \cite[Propositions 3.3 and 3.5]{moakherMeansAveragingGroup2002} that $\mu_{F}(\tilde {\mf g}; \mathcal D)$ can also be viewed as a projected mean and is readily computed as $\mc P_{SO}(\mu_E(\tilde {\mf R}))$, where $\mc P_{SO}$ is defined identically to the $SO(2)$ case. Therefore, the Fr\'echet mean is almost surely (i.e., except in the measure-zero event that $\mu_E(\tilde {\mf R})$ is singular) unique and computable in closed-form.

As explained in \Cref{sec:relationship}, the Frobenius inner product $\ip{\mf X}{\mf Y}_{\mathfrak{so}(3)}\coloneqq\tr(\mf X \tran\mf Y)$ can be uniquely extended to define a bi-invariant Riemannian metric $\mathcal R$ on $SO(3)$, whose geodesic distance function is $\mc R(\mf R,\mf Q)=\lVert \log(\mf R\tran\mf Q)\rVert$.
Any point $\bmu\in SO(3)$ that satisfies
$$
\int_{SO(3)}\log(\bmu\tran\mf R) f(\mf R) d\mf R=\mf 0_3
$$
is an element of $\mu_G(\tilde{\mf R})$. Moreover, since $\mu_G(\tilde{\mf R})=\mu_K(\tilde{\mf R}; \mc R)$, it follows that $\mu_F(\tilde {\mf R}; \mathcal R)\subseteq \mu_G(\tilde{\mf R})$. 




\subsection{$SE(2)$ and $SE(3)$}
\label{sec:SEn}
In this section, we consider the \textit{Special Euclidean} group $SE(d)$, focusing on the cases of $d=2$ and $d=3$.
As a manifold, $SE(d)$ is diffeomorphic to the direct product group, $SO(d)\times \mathbb R^d$. Each element in $SO(d)\times \mathbb R^d$ is of the form $(\mf R, \mf t)$ (i.e., a pair of rotation matrix and translation vector). The map $\mf H$ introduced in \Cref{sec:extrinsic} is a smooth, invertible map from $SO(d)\times \mathbb R^d$ to $SE(d)$. As a Lie group, $SE(d)$ is a \textit{semi-direct product} of the constituent groups (written as $SE(d)=SO(d)\ltimes \mathbb R^d$), and has a different group structure from the direct product. The group operation of $SO(d)\times \mathbb R^d$ is given by
\begin{align}
(\mf R, \mf t)(\mf Q, \mf s) = (\mf R \mf Q,\, \mf t + \mf s),
\end{align}
which is different from the group operation of $SE(d)$ that was illustrated in \cref{eq:se3-multiplication}.
Given an element $\mf g\in SE(3)$, we will write $\mf R_{\mf g}$ and $\mf t_{\mf g}$ to refer to its rotation and translation parts, so that $\mf H\left(\mf R_{\mf g}, \mf t_{\mf g}\right)= \mf g$.

Let $\tilde {\mf g}$ be a random variable on $SE(d)$. Observe that $\mf R_{\widetilde {\mf g}}$ is a random variable on $SO(d)$ whereas $\mf t_{\widetilde {\mf g}}$ is a random variable on $\mathbb R^d$. 
Hence, the Euclidean mean of $\tilde{\mf g}$ is a function of the Euclidean means of $\mf R_{\widetilde {\mf g}}$ and $\mf t_{\widetilde {\mf g}}$:
\begin{align}
\bmu_E(\tilde{\mf g}) &=
\int_{SE(d)} 
\begin{bmatrix}
    \mf R_{\mf g} & \mf t_{\mf g}\\
    \mf 0_{1\times d} & 1
\end{bmatrix} f(\mf g)d\mf g
=
\begin{bmatrix}
   \bmu_E(\mf R_{\widetilde {\mf g}}) &  \bmu_E(\mf t_{\widetilde {\mf g}}) \\
\mf 0_{1\times d} & 1
\end{bmatrix}.
\end{align} 
Define ${\mf g} \cdot\mf x \coloneqq \mf R_{\mf g} \mf x + \mf t_{\mf g}$, where $\mf g\in SE(d)$ and $\mf x\in \mathbb R^d$. This makes $\mathbb R^d$ a \textit{homogeneous space} equipped with a left-action of $SE(d)$, such that the group element $\mf g$ \textit{acts} on $\mf x$ to yield $\mf g\cdot \mf x$.
As discussed in \Cref{sec:extrinsic}, 
$\bmu_E\big(\tilde{\mf g})\cdot\mf x$ is the mean of the random vector $\tilde{\mf g} \cdot\mf x$. This means that the Euclidean mean is a suitable definition of mean when one is interested in the action of $\tilde{\mf g}$ on the vectors of $\mathbb R^d$, viewed as a homogeneous space. For instance, $\tilde{\mf g}$ can describe the uncertain pose of a rigid body, and $\mf x$ may be chosen as the centroid of the corresponding body. In this case, the random vector $\tilde{\mf g}\cdot \mf x$ describes the uncertainity in the centroid of the rigid body, so that $\bmu_E(\tilde{\mf g}\cdot \mf x)=\bmu_E(\tilde{\mf g})\cdot \mf x$ is indeed a quantity of interest. 

\begin{remark}
The adjoint representation of $SE(3)$ can be viewed as a matrix Lie group. Specifically, we define $\bm\Ad_{SE(3)}\coloneqq\lbrace\bm\Ad_{\mf g}|\mf g\in SE(3)\rbrace$ (with the basis for $\mathfrak{se}(3)$ chosen as in \cite{chirikjian2011stochastic}), which is isomorphic to $SE(3)$ \cite{barfoot2024state}. 
However, $\bm\Ad_{\bmu_E(\tilde{\mf g})}\neq \bmu_E(\bm\Ad_{\tilde{\mf g}})$ in general, showing that the Euclidean mean explicitly depends on the choice of representation.
\end{remark}

Since $\bmu_E(\mf R_{\widetilde {\mf g}}) \notin SO(d)$, the Euclidean mean of $\tilde{\mf g}$ is not a homogeneous transformation matrix. Nevertheless, one can use the function $\mc P_{SO}$ to project the mean back to $SE(d)$. That is, the transformation matrix $\mf H\left(\mf Q, \mf s\right)$, where $\mf Q\coloneqq \mc P_{SO}(\bmu_E(\mf R_{\widetilde {\mf g}}))$ and $\mf s\coloneqq \bmu_E(\mf t_{\widetilde {\mf g}})$, is a projected mean of $\tilde{\mf g}$. The projected mean has the advantage of being computable in closed-form, and can serve as a good initial guess for the purpose of computing the group-theoretic, Karcher, or Riemannian Fr\'echet means using iterative algorithms. Recently, some authors have used \textit{dual quaternions} (which are an $8$-dimensional extension of quaternions) to represent random variables on $SE(3)$ \cite{srivatsan2016estimating,li2020unscented}. However, the averaging of dual quaternions suffers from drawbacks similar to those of quaternion-based averaging on $SO(3)$ when $\tilde{\mf g}$ has a large variance.

The group-theoretic means on $SE(2)$ and $SE(3)$ have been used for multi-robot localization \cite{li2016lie}, analysis of robotic manipulators \cite{wang2006error}, and multiscale modeling of DNA \cite{wolfe2012multiscale}. They were investigated from the standpoints of existence and uniqueness in \cite{pennec2013bi,pennec2020beyond}.
However, since these groups do not admit a bi-invariant metric, the group-theoretic means do not minimize the cost function given in \cref{def:lie-karcher-mean}. Nevertheless, some authors minimize cost functions similar to \cref{def:lie-karcher-mean} to obtain a notion of mean (that is neither the group-theoretic mean nor a Fr\'echet mean) for a random variable in $SE(2)$ or $SE(3)$ \cite{barfoot2025certifiably}.

In the absence of a bi-invariant metric, a natural choice for defining the Riemannian Fr\'echet and Karcher means is the left-invariant metric induced by a choice of inner product on $\mathfrak {se}(d)$. 
As shown in \cite{park1995distance}, if one uses the Frobenius inner product $\ip{\mf X}{\mf Y}\coloneqq \tr(\mf X\tran \mf Y)$ on $\frak{se}(3)$, then the corresponding geodesic distance function is given by
\begin{align}
\mc R(\mf g, \mf h)&=\sqrt{\lVert\log\big({\mf R_{\mf g\,}\tran} \mf R_{\mf h}\big)\rVert^2 + \lVert \mf t_{\mf g} - \mf t_{\mf h}\rVert^2},
\end{align}
wherein $\log(\,\cdot\,)$ is the $\log$ map of $SO(3)$. Interestingly, this is also the geodesic distance function of $SO(3)\times \mathbb R^3$ (when endowed with the product of the bi-invariant metrics on $SO(3)$ and $\mathbb R^3$). 

In addition to the Riemannian distance, one can define Fr\'echet means on $SE(d)$ using more general distance functions that might be better suited for certain applications. For instance, one can define:
\begin{align}
\mc D_{B}(\mf g, \mf h) &= \sqrt{\lVert \mf R_{\mf g} - \mf R_{\mf h} \rVert^2 + m \lVert \mf t_{\mf g} - \mf t_{\mf h} \rVert^2}.
\end{align}
As shown in 
\cite[p. 208]{chirikjian2016harmonic}, such a distance function arises naturally in the case where $\tilde {\mf g}$ describes the uncertain pose of a rigid body.
In this case, the Fr\'echet mean
is also right-invariant, and has the advantage of being unique and computable in closed-form using the polar decomposition $\mc P_{SO}$. That is, the set of Fr\'echet means $\mu_F(\tilde{\mf g};\mc D_B)$ is (almost surely) the singleton set  $\lbrace\mf H(\mf Q, \mf s)\rbrace$, where $\mf Q\coloneqq \mc P_{SO}(\bmu_E(\mf R_{\widetilde {\mf g}}))$ and $\mf s\coloneqq\bmu_E(\mf t_{\widetilde {\mf g}})$.
It also offers a way to interpolate between $SE(3)$ elements, yieding curves that are not geodesics, but are nonetheless of practical relevance \cite{liu2024average}.

\section{Conclusion}\ 
In this article, we investigated the problem of defining notions of \textit{mean} and \textit{covariance} for a random variable that takes values in a matrix Lie group. 
The existing definitions of means on Lie groups were reviewed, and their relationships to each other were clarified. We also discussed the various properties of these means, including their dependence on additional geometric structure on the Lie algebra (such as an inner product and/or choice of basis).
The \textit{Special Orthogonal} and \textit{Special Euclidean} Lie groups were considered as motivating examples, as these Lie groups are frequently encountered in engineering applications.
By highlighting the subtleties involved in defining and computing means on noncommutative Lie groups, this article lays the foundation for principled statistical inference on Lie groups and serves as a guide for practitioners faced with the task of selecting an appropriate definition of mean for a given application.
\appendix
\section{Review of Lie Groups} 
\label{app:mlg}

\subsection{Matrix Lie Groups}
A \textit{matrix Lie group} refers to an $n$-dimensional Lie group equipped with an $m$-dimensional matrix representation. For instance, ${GL}(m, \mathbb C)$ is the group of all $m\times m$ invertible matrices with complex entries, and has the dimension $n=m^2$. The Lie algebra of ${GL}(m,\mathbb C)$, denoted as $\mathfrak {gl}(m,\mathbb C)$, is the space $\mathbb C^{m\times m}$ (viewed as a vector space under addition) equipped with the Lie bracket $[\mf A, \mf B]\coloneqq \mf A \mf B - \mf B \mf A$.
We write $G \leq {GL}(m, \mathbb C)$ to indicate that $G$ is a Lie subgroup of ${GL}(m, \mathbb C)$. It can be shown that $\mathfrak g$ is a subset of $\mathfrak {gl}(m,\mathbb C)$ that is closed under the Lie bracket. We say that $\mathfrak g$ is \textit{real} if it is real as a vector space; equivalently, if it is a subset of $\mathfrak {gl}(m,\mathbb C)$ that is closed under vector addition, scalar multiplication by real numbers, and the Lie bracket \cite{hall2013lie}. 


The \textit{exponential map} $\exp:\mathfrak g \rightarrow G$ relates each vector $\mf X\in\mathfrak g$ to a group element $\exp(\mf X)\in G$. Since we have assumed that $G$ is a matrix Lie group, $\exp$ coincides with the matrix exponential 
$\exp(\mf X) = \sum_{i=1}^\infty \frac{1}{i!}\mf X^i$ whenever this sum converges. There exists an open neighborhood of $\mf 0_m\in \mathfrak g$, denoted as $\subs{\mathfrak g}\subseteq \mathfrak g$, where the exponential map has a well-defined inverse, which we call the (Lie-theoretic) {log map} of $G$ \cite{gallier2020differential,chirikjian2011stochastic}, also known as the \text{principal logarithm} \cite{hall2013lie}. 
The set $\subs{G}\coloneqq \exp(\subs{\mathfrak g})$ represents the domain of the {log map};
for $\mf X\in \subs{\mathfrak g}$, it holds that $\log\exp(\mf X) = \mf X$, and for $\mf g\in \subs{G}$, it holds that $\exp\log (\mf g) = \mf g$.
The $\exp$ and $\log$ maps (in general) are not the same as the Riemannian $\rexp{\mf g}$ and $\rlog{\mf g}$ maps introduced in Appendix \ref{app:riemannian}.
The exponential and log maps for some Lie groups (including $SO(d)$ and $SE(d)$) can be computed using closed-form formulae \cite{chirikjian2011stochastic,barfoot2024state}.

Let $G$ be a real $n$-dimensional Lie group, and $\lbrace \mf E_i\rbrace_{i=1}^n$ a basis for $\mathfrak g$. Given a vector $\mf X\in \mathfrak g$, it can be uniquely expressed as $\mf X=X^i \mf E_i$, where $\mf X^i\in \mathbb R$.\footnote{Here (and in the rest of the paper, unless indicated otherwise), the Einstein summation convention is used, so that $\mf X^i \mf E_i=X^1\mf E_1 + X^2 \mf E_2+\ldots+X^n \mf E_n$.}
Hence, the choice of basis defines the map $(\,\cdot\,)^\vee: \mathfrak g \rightarrow \mathbb R^n$ and its inverse $(\,\cdot\,)^\wedge$ (read as `vee' and `hat', respectively), which identify $\mathfrak g$ with the vector space $\mathbb R^n$:
\begin{align}\label{def:vee}
(X^i \mf E_i)^\vee = \begin{bmatrix}
X^1\\\vdots\\X^n
\end{bmatrix} \qquad \text{and} \qquad \left(\begin{bmatrix}
X^1\\\vdots\\X^n
\end{bmatrix}\right)^\wedge = X^i \mf E_i.
\end{align}
The \textit{structure constants} of $\mathfrak g$ in the basis $\lbrace \mf E_i\rbrace_{i=1}^n$ are denoted as $\lbrace C_{ij}^k\rbrace_{i,j,k=1}^n$, and are defined by the equality $[ \mf E_i,  \mf E_j]= C_{ij}^k  \mf E_k$. The \textit{adjoint representation} of $G$ maps each element to an invertible linear transformation of $\mathfrak g$; given $\mf g\in G$, we have $\Ad_{\mf g} \mf X = \mf g \mf X \mf g^{-1}$.
In the basis $\lbrace \mf E_i \rbrace_{i=1}^n$, we can express it as the matrix $\bm \Ad_{\mf g} \in \mathbb R^{n\times n}$, defined such that $\left(\Ad_{\mf g} \mf X\right)^\vee = \bm \Ad_{\mf g} \mf X^\vee$. Relatedly, there is an adjoint representation of $\mathfrak g$, given by
$$ \ad_{\mf X} \mf Y\coloneqq \frac{d}{dt}\Ad_{\,\exp(t \mf X)}\mf Y\bigg|_{t=0}=[\mf X,\mf Y].$$
The matrix of $\ad_{\mf X}$ is $\bm \ad_{\mf X}$, so that $(\ad_{\mf X}\mf Y)^\vee =\bm\ad_{\mf X}\mf Y^\vee$.
Note that the columns of $\bm \ad_{\mf X} \in \mathbb R^{n\times n}$ are defined via the action of $\ad_{\mf X}$ on each of the basis vectors:
\begin{align}\label{eq:ad-matrix-columns}
\bm\ad_{\mf X}&\coloneqq\begin{bmatrix}
(\ad_{\mf X} \mf E_1)^\vee& (\ad_{\mf X} \mf E_2)^\vee& \cdots & (\ad_{\mf X} \mf E_n)^\vee
\end{bmatrix}.
\end{align}
Using (\ref{eq:ad-matrix-columns}), it can be shown that the $(i,j)^{th}$ component of $\bm\ad_{\mf X}$ is $(\bm\ad_{\mf X})_{ij}=X^k  C_{kj}^i.$
The adjoint representations respect the group multiplication of $G$ and the Lie bracket of $\mathfrak g$:
\begin{align}
\bm\Ad_{\mf g \mf h}&=\bm\Ad_{\mf g} \bm\Ad_{\mf h}\quad\text{and}\quad\bm\ad_{[\mf X,\mf Y]} = [\bm\ad_{\mf X}, \bm\ad_{\mf Y}]
\end{align}
for all $\mf g,\mf h\in G$ and $\mf X,\mf Y\in\mathfrak g$.

\subsection{Inner Products on the Lie Algebra}
\label{app:inner-products}
Let $\lbrace \mf E_i \rbrace_{i=1}^n$ be a basis  and $\ip{\,\cdot\,}{\,\cdot\,}_{\frak g}$ an arbitrary \textit{inner product} for $\mathfrak g$.
The {inner product} can be described using an $n\times n$ positive definite matrix $\mf W$, whose components are defined as
\begin{align}
\mf W_{ij}&\coloneqq \ip{\mf E_i}{\mf E_j}_{\frak g}\quad\forall i,j\in\lbrace1,\cdots,n\rbrace.\label{def:weighted-inner-product}
\end{align}
Let $\sqrt{\mf W}$ be the symmetric square root of $\mf W$ (whose existence and uniqueness follows from \cite[Thm. 7.2.6]{hornMatrixAnalysis2012}), i.e., $\sqrt{\mf W}^2 = \mf W$ and $\sqrt{\mf W}\tran = \sqrt{\mf W}$. Then, we can define the basis $\lbrace {\mf E}_i^\prime \rbrace_{i=1}^n$ via 
$\mf E_i^\prime \coloneqq \sqrt{\mf W}^{ij}\mf E_j,$ 
where $\sqrt{\mf W}^{ij}\coloneqq \big(\sqrt{\mf W}^{-1}\big)_{ij}$. It is straightforward to verify that $\lbrace \mf E_i^\prime \rbrace_{i=1}^n$ is orthonormal w.r.t.~the foregoing inner product, i.e., $\ip{\mf E_i^\prime}{\mf E_j^\prime}_{\frak g}=\delta_{ij}$. Therefore, no loss of generality is caused by choosing to work with an orthonormal basis.
However, one often knows the structure constants and their properties in the original basis $\lbrace \mf E_i \rbrace_{i=1}^n$, in which they take a simpler form than they do in the orthonormal basis. For instance, the structure constants of $\mathfrak{so}(3)$, in the choice of basis for $\mathfrak{so}(3)$ considered in \Cref{sec:SU2}, are the well-known permutation symbols (also known as Levi-Civita symbols) from multivariable calculus \cite[Prop. C7]{hall2013lie}. Therefore, it is worthwhile to discuss how various quantities (e.g., structure constants) transform under the change of basis, $\lbrace \mf E_i \rbrace_{i=1}^n\mapsto \lbrace \mf E_i^\prime \rbrace_{i=1}^n$.

Given a vector $\mf X \in \mathfrak g$, it can be expressed in either basis as $\mf X=X^i \mf E_i={X^\prime}^i \mf E_i^\prime$. 
To the basis $\lbrace \mf E'_i \rbrace_{i=1}^n$ is a corresponding ${(\,\cdot\,)^\vee}'$ map, defined analogously to (\ref{def:vee}). 
Similarly, let $C'^k_{ij}$ denote the structure constants in the $\lbrace \mf E'_i \rbrace_{i=1}^n$ basis. The coefficients in either basis are related via the following formulae:
\begin{align}
{X^\prime}^i = \sqrt{\mf W}_{ij}{X}^j\quad \text{and} \quad
 {C'^k_{ij}} = \sqrt{\mf W}^{is}\sqrt{\mf W}^{j\ell} \sqrt{\mf W}_{kr} C^r_{s\ell}. \label{eq:struct-const-transform}
%
\end{align}
Let $\bm \ad'_{\mf X}$ be the matrix of $\ad_{\mf X}$ in the $\lbrace E'_i \rbrace_{i=1}^n$ basis. Using the expression for $C'^i_{kj}$, we have
\begin{align}
({\bm\ad}_{\mf X})^\prime_{ij}={X^\prime}^kC'^i_{kj} 
=\sqrt{\mf W}^{j\ell} \sqrt{\mf W}_{ir} C^r_{s\ell} = {X}^s\sqrt{\mf W}^{j\ell} \sqrt{\mf W}_{ir} (\bm\ad_{\mf X})_{r\ell}.
\label{eq:ad-components-transform}
\end{align}
Equation (\ref{eq:ad-components-transform}) can be rewritten as $\bm\ad'_{\mf X}=\sqrt{\mf W}\bm\ad_{\mf X}\sqrt{\mf W}^{-1}$, which shows that $\bm\ad_{\mf X}$ undergoes a similarity transformation under the change of basis. In fact, it can be shown that the matrix of \textit{any} linear transformation from $\mathfrak g$ to $\mathfrak g$ changes by a similarity transformation under a change of basis, so we have $\bm\Ad'_{\mf g} = \sqrt{\mf W}\bm\Ad_{\mf X}\sqrt{\mf W}^{-1}$ as well. 

\section{Calculus on Lie groups}
\label{app:mlg-calc}
Given a vector $\mf Z\in \mathfrak g$ and a group element $\mf g\in G$, consider the curve $\bm\upgamma:\mathbb R \rightarrow G$ defined by $\bm\upgamma(t)=\mf g\exp(t\mf Z)$. At $t=0$, the curve $\bm\upgamma(t)$ has the value $\bm\upgamma(0) = \mf g$ and the derivative (i.e., tangent vector) $\dot {\bm\upgamma}(0) = \mf g\mf Z$. The space $T_{\mf g}G\coloneqq \lbrace \mf g\mf Z\,|\,\mf Z\in \frak g\rbrace$ can be given a vector space structure, at which point it is called the tangent space of $G$ at $\mf g$. Note that $T_{\mf I_m} G$ is the same as $\mathfrak g$ (except that $\mathfrak g$ has, in addition, the Lie bracket operation).

\subsection{Parametrization}\label{app:parametrization}
A \textit{local parametrization} of $G$ refers to a diffeomorphism (i.e., a smooth invertible map) $\psi:\subs{\mathbb R^n} \rightarrow S_{G}$, where $\subs{\mathbb R^n}\subseteq \mathbb R^n$ and $S_{G}\subseteq G$, $n$ being the dimension of $G$ as a manifold. 
The map $\psi^{-1}$ is then a system of coordinates (i.e., a chart) for $S_{G}$. 
For example, the parametrization for $SO(2)$ described in \cref{eq:so2-param} is defined as the map $\mf R:(-\pi,\pi)\rightarrow SO(2)\backslash\lbrace-\mf I_2\rbrace$. Letting $\mf x\in \subs{\mathbb R^n}$ be a point in the parameter space,
we define the \textit{Jacobian} of $\psi$ as the matrix
\begin{align}
     \mf J_{\psi}(\mf x) \coloneqq \begin{bmatrix}\left(\psi(\mf x)^{-1}\frac{\partial \psi}{\partial x^1}(\mf x)\right)^\vee & \left(\psi(\mf x)^{-1}\frac{\partial \psi}{\partial x^2}(\mf x)\right)^\vee & \cdots & \left(\psi(\mf x)^{-1}\frac{\partial \psi}{\partial x^n}(\mf x)\right)^\vee\end{bmatrix}.
\end{align}
This matrix is the differential (i.e., the pushforward map) of $\psi$ at ${\mf x}$, as expressed w.r.t. the standard basis for $T_{\mf x}\mathbb R^n$ and the basis $\lbrace {\psi(\mf x)\mf E_{i}}\rbrace_{i=1}^n$ for $T_{\psi(\mf x)}G$. 

The Jacobian of the inverse map $\psi^{-1}$ is given by $\mf J_{\psi^{-1}}(\mf g)=\mf J_{\psi}\big(\psi^{-1}{\small (\mf g)}\big)^{-1}$. As in the case of the $\exp$ and $\log$ maps, the matrices $\mf J_{\psi}(\mf x)$ and $\mf J_{\psi^{-1}}(\mf g)$ can be computed using closed-form expressions in the cases of $G=SO(d)$ or $SE(d)$ \cite{barfoot2024state}.


\subsection{Differentiation}
We can extend $\mf Z$ to a unique \textit{left-invariant vector field (LIVF)} $\ms Z$ on $G$, which assigns to each point $\mf g\in G$ the tangent vector
$\ms Z_{\mf g} \coloneqq\mf g \mf Z$. It has the property $\ms Z_{\mf h\mf g} = \mf h \ms Z_{\mf g}$ for all $\mf h\in G$, which explains the term `left-invariant' (i.e., invariance under left-multiplication). Moreover, it holds that $\dot{\bm\upgamma}(t)=\ms Z_{\bm\upgamma(t)}$, due to which $\bm\upgamma$ is called the \textit{integral curve} of $\ms Z$ starting at $\mf g$. Given a smooth function $f:G\rightarrow V$ where $V$ is a vector space, the \textit{Lie derivative} of $f$ along $\ms Z$ is defined as
\begin{align}
\big(\ms Z f\big)(\mf g) &\coloneqq \frac{d}{dt}f\big(\mf g\exp(t\mf Z)\big)\Big\vert_{t=0}= \lim_{\tau \rightarrow 0}\frac{f\big(\mf g\exp(\tau \mf Z)\big) - f(\mf g)}{\tau}.
\label{def:lie-derivative}
\end{align}
In this sense, $\ms Z$ is a \textit{differential operator} that maps $f$ to another function on $G$; the latter function represents the derivative of $f$ along $\bm\upgamma$. 
Similarly, one can define a right-invariant vector field (RIVF) and a corresponding differential operator. We remark that $\ms Z$ is the same as the \textit{right} Lie derivative in \cite{chirikjian2011stochastic} (where it is denoted as $\mf Z^{(r)}$), with `right' referring to the fact that $\exp(\,\cdot\,)$ appears to the {right} of $\mf g$ in (\ref{def:lie-derivative}).

\subsection{Integration}
The concept of integration on $G$ can be defined w.r.t. the \textit{left Haar measure}. It is characterized (up to scaling) by the following \textit{left-invariance} property.
If $S\subseteq G$ and $\mf hS \coloneqq \lbrace \mf h \mf g \,|\, \mf g\in S\rbrace$, then
\begin{align}\label{def:left-haar}
\int_S d \mf g &= \int_{\mf hS} d \mf g \quad \text{for all } \mf h\in G,
\end{align}
where `$d\mf g$' denotes integration w.r.t.~the (left) Haar measure.
If $G$ is compact, then the integral in \cref{def:left-haar} is finite, and there exists a unique \textit{normalized} Haar measure for which the integral $\int_G d\mf g$ is equal to $1$.

The volume of $S\mf h$ is not the same as that of $S$ in general; when this is true, the Lie group is said to be \textit{unimodular}, and we say that the left Haar measure is \textit{bi-invariant}. An important consequence of bi-invariance of the Haar measure is that given 
any (measurable) function $f:G\rightarrow V$, where $V$ is a vector space, we have
\begin{align}
\int_G f(\mf g) d \mf g = \int_G f(\mf h \mf g) d \mf g = \int_G f(\mf g\mf h) d \mf g = \int_G f(\mf g^{-1}) d \mf g,
\end{align}
for all $\mf h\in G$. The condition for $G$ to be \textit{unimodular} is $\det (\Ad_{\mf g})= 1$ for all $\mf g\in G$. Writing $\mf g$ as $\exp(\mf X)$ for some $\mf X\in \mathfrak g$, we have that \cite[Ch. 4]{faraut2008analysis}\begin{align}\label{eq:unimodularity-Ad-ad}
\det (\Ad_{\exp(\mf X)})=\det \big(\exp(\ad_\mf X)\big)=e^{\tr(\ad_\mf X)}=1.\end{align}
Thus, $G$ is unimodular if and only if $\tr(\ad_\mf X)=0$ for all $\mf X\in \mathfrak g$. The conditions in \cref{eq:unimodularity-Ad-ad} can be checked by choosing an arbitrary basis for $\mathfrak g$, since the determinant and trace of a linear operator are both basis-independent concepts.
In terms of the structure constants (w.r.t. an arbitrary basis for $\mathfrak g$), the condition for unimodularity is that $C_{ij}^j=0$ for all $i=1,\ldots,n$.

To compute the Haar integral numerically, one typically uses the parametrization $\psi:\subs{\mathbb R^n}\rightarrow \subs{G}$, where we recall that $\subs{\mathbb R^n}\subseteq \mathbb R^n$ and $\subs{G}\subseteq G$. 
The function $f$ and the Haar measure must each be \textit{pulled back} under $\psi$ to make the integrals on $\subs{\mathbb R^n}$ and $\subs{G}$ coincide, giving us \cite{faraut2008analysis,chirikjian2011stochastic}
\begin{align}\label{eq:integrals-parametrization}
\int_{\subs{G}}f(\mf g)d \mf g &= \int_{\subs{\mathbb R^n}}f\left(\psi(\mf x)\right) \,\vert \det\left(\mf J_{\psi}(\mf x)\right)\vert\,d\mf x,
\end{align}
where $|\,\cdot\,|$ represents the absolute value of a real number. 
Examples of Haar integrals for commonly-encountered Lie groups can be found in \cite{chirikjian2011stochastic}.

\subsection{Probability Density Functions}
In this paper, we assume that the probability density function (pdf) of a random variable $\tilde {\mf g}$, is defined w.r.t the Haar measure. Letting $f:G \rightarrow \mathbb R_{\geq 0}$ be such a pdf, we have $\int_G f(\mf g) d \mf g = 1$. If the support of $f$ is contained in $\subs{G}$, then it can be pulled back to define the pdf $$\bar f(\mf x) \coloneqq  f(\psi(\mf x))\,\big\vert \det\left(\mf J_{\psi}(\mf x)\right) \big\vert$$
on the parameter space $\subs{\mathbb R^n}$. Using \cref{eq:integrals-parametrization}, we see that the pdf $\bar f$ can be integrated on $\subs{\mathbb R^n}$ to get $1$, as expected. 
Conversely, if $\bar f(\mf x)$ is the pdf of a random variable $\tilde{\mf x}$ in $\subs{\mathbb R^n}$, then $$f(\mf g)\coloneqq \bar f\left(\psi^{-1}(\mf g)\right)\vert\det(\mf J_{\psi^{-1}}(\mf g))\vert$$ is the corresponding pdf on $\subs{G}$ that makes the equality in \cref{eq:integrals-parametrization}
hold.%
\footnote{This can be shown by using the fact that $\det(\,\mf J_{\psi^{-1}}(\mf g))
={\det(\mf J_{\psi}(\mf x))}^{-1},$ where $\mf g=\psi(\mf x).$}
As random variables can be characterized (up to measure-zero events) by their pdfs, random variables in $\subs{\mathbb R^n}$ are equivalent to random variables on $\subs{G}$ and vice versa.

Let $C^\infty(G)$ be the space of smooth functions on $G$. The \textit{Dirac delta function} $\delta(\mf g)$ can be rigorously defined as a \textit{distribution}, i.e., an element of the vector space dual of $C^\infty(G)$ \cite[Sec. 2.18]{schutz1980geometrical} \cite[p.298]{varadarajan1999introduction}. For the purposes of this paper, it is sufficient to characterize $\delta(\mf g)$ via the property
\begin{align}
    \int_{G} f(\mf g)\delta(\mf g) d \mf g = f(\mf I_m) \quad \forall\,f\in C^\infty(G),
\end{align}
which is consistent with the use of Dirac delta functions in the physics literature.
Given a set of $N$ samples on $G$, $\lbrace \mf g_i \rbrace_{i=1}^N$, one can describe this set using the \textit{empirical pdf} $f(\mf g)=\frac{1}{N}\sum_{i=1}^N\delta(\mf g_i^{-1}\mf g)$. The Euclidean mean of $\lbrace \mf g_i \rbrace_{i=1}^N$ (as defined in Section \ref{sec:euclidean}) can be computed using the empirical pdf:
\begin{align}
\int_G \mf g \bigg(  
\frac{1}{N}\sum_{i=1}^N\delta(\mf g_i^{-1}\mf g)
\bigg)
d \mf g &= 
\frac{1}{N}\sum_{i=1}^N
\bigg(\int_G \mf g\,   
\delta(\mf g_i^{-1}\mf g)
d\mf g\bigg) \\
&=
\frac{1}{N}\sum_{i=1}^N
\bigg(\int_G \,\mf g_i \mf g\,   
\delta(\mf g)
d\mf g\bigg)
=
\frac{1}{N}\sum_{i=1}^N \mf g_i,
\end{align}
where we have used the linearity and left-invariance of the Haar integral.
In this way, Dirac delta functions enable us to treat continuous pdfs and empirical pdfs on $G$ in an equal footing, without introducing advanced measure-theoretic concepts.

\section{Riemannian Geometry on Lie Groups}
\label{app:riemannian}
The inner product $\langle \cdot , \cdot \rangle$ on $\mathfrak{g}$ induces a unique \textit{left-invariant Riemannian metric} $\langle \cdot , \cdot \rangle^L_{\mf g}$ on $G$, defined via the relation
\begin{align}\label{def:left-invariant-metric}
\langle \mf V, \mf W \rangle^L_{\mf g} \coloneqq \langle  \mf g^{-1}\mf V,  \mf g^{-1}\mf W \rangle
\end{align}
\noindent 
for all $\mf V, \mf W\in T_{\mf g} G$. Roughly speaking, $\langle \cdot , \cdot \rangle^L_{\mf g}$ is an inner product for the tangent space $T_{\mf g} G$ that is also \textit{smooth} as a function of $\mf g$. If one inserts a pair of smooth vector fields into $\langle \cdot, \cdot \rangle^L$, then the result is a smooth function on $G$. The term \textit{left-invariance} refers to the property $\langle \mf V, \mf W \rangle^L_{\mf g}=\langle \mf h\mf V, \mf h\mf W \rangle^L_{\mf h\mf g}$. 

Each choice of Riemannian metric on $G$ offers a way to differentiate tangent vectors, via the \textit{Levi-Civita connection} \cite{lee2012smooth,cheegerComparisonTheoremsRiemannian1975a}. One says that a curve on $G$ is a \textit{geodesic} if the derivative of its velocity (as computed using the Levi-Civita connection) is zero. In the special case of $G=(\mathbb R^n, +)$, the geodesics of a left-invariant Riemannian metric are straight lines, and vice versa. An alternative characterization of geodesics on Lie groups is given in \Cref{app:geodesics}, which circumvents the use of a connection.

\subsection{Bi-Invariant Riemannian Metrics}
We say that $\langle \cdot , \cdot \rangle^L_{\mf g}$ is \textit{right-invariant} if it holds that
$\langle \mf V, \mf W \rangle^L_{\mf g} = \langle \mf V \mf h, \mf W\mf h \rangle_{\mf g\mf h}^L$,
in which case it is said to be \textit{bi-invariant}. 
The following lemma lists the conditions under which the left-invariant metric $\langle \cdot , \cdot \rangle^L_{\mf g}$ is bi-invariant.

\begin{lemma}\label{lem:bi-invariant-metrics}
Let an inner product on $\mathfrak g$ be extended to define a left-invariant Riemannian metric on $G$. Assume that $G$ is connected. Then, the following are equivalent:
\begin{enumerate}[\textup{(\alph*)}]
    \item The inner product is $\Ad$-invariant, i.e.,  for all  $\mf X,\mf Y\in \mathfrak g$  and  $\mf g\in G$,
    $$\langle \mf X, \mf Y \rangle = \langle \Ad_{\mf g} \mf X,\,  \Ad_{\mf g} \mf Y \rangle. $$
    \item The map $\ad_\mf X$ is a skew-isometry, i.e.,  for all $\mf X,\mf Y,\mf Z\in \mathfrak g$,
    $$\langle \mf Y, \ad_\mf X \mf Z \rangle = -\langle \ad_\mf X \mf Y, \mf Z \rangle.$$
    \item The left-invariant Riemannian metric is also right-invariant, i.e., it is bi-invariant.
    \item Any curve of the form $\bm\upgamma(t) = \mf g\exp(t \mf X)$, where $\mf g\in G$ and $\mf X\in\frak g$, is a geodesic. 
    \item It holds that $\bm\ad_{\mf x^\wedge}\tran \mf W\mf x=\mf 0_{n\times 1}$ for all $\mf x\in\mathbb R^n$.
    \item The matrix ${\mf C_\ell}$ whose $(i,j)^{th}$ component is defined by $({\mf C_\ell})_{ij} = C_{i\ell}^k \mf W_{kj}$ is skew-symmetric.%
\end{enumerate}
\end{lemma}
\begin{proof}
It was shown by Milnor in \cite{milnorCurvaturesLeftInvariant1976} that $\Ad$-invariance of the inner product, $\ad_\mf X$ being a skew-isometry, and bi-invariance of the corresponding left-invariant metric are equivalent conditions; this is \textup{(a)}, \textup{(b)}, and \textup{(c)}, respectively.
The equivalence of \textup{(d)} and \textup{(e)} is the content of \cite[Prop. 3.18]{cheegerComparisonTheoremsRiemannian1975a}. The remaining implications will be shown here.

First, we show that $\textup{(b)} \Rightarrow \textup{(e)}$. The \textit{adjoint} of the operator $\ad_\mf X$ (with respect to the inner product on $\mathfrak g$) is the unique operator $\ad^\ast_\mf X$ that satisfies $\ip{\mf Y}{\ad_\mf X \mf Z}=\ip{\ad^\ast_\mf X \mf Y}{\mf Z}$. 
It is readily shown that the matrix of $\ad_{\mf X}^\ast$ satisfies $\mf W\bm\ad_\mf X = (\bm\ad_{\mf X}^\ast)\tran \mf W$. Since $\bm\ad_{\mf X}^\ast = -\bm\ad_\mf X$ (using $(\textup b)$), we have $\mf W\bm\ad_\mf X = -\bm\ad_\mf X\tran \mf W$. Condition $(\textup e)$ then follows from the fact that $\ad_\mf X \mf X=[\mf X,\mf X]=\mf 0_m$. 

Next, we show that $\textup{(e)}\Rightarrow \textup{(f)}$. Using the linearity of the adjoint representation, we can rewrite \textup{(e)} as $
X^i (\bm\ad_{\mf E_i}\tran)_{\ell j} \mf W_{jk} X^k =0
$
$\forall \ell \in \lbrace 1, 2, \ldots, n\rbrace$.
In terms of the structure constants, we have
\begin{align}
X^i C_{i \ell}^j \mf W_{jk} X^k =0 \ \Rightarrow\ \mf x^\top \mf C_{\ell}\mf x = 0.
\end{align}
It follows that $\mf x\tran  ({\mf C_\ell} + {\mf C_\ell}\tran )\,\mf x = 0$ for all $\mf x \in\mathbb R^n$. Since $({\mf C_\ell} + {\mf C_\ell}\tran )$ is symmetric, the preceding condition reduces to ${\mf C_\ell} + {\mf C_\ell}\tran  = \mf 0$, i.e., ${\mf C_\ell}$ is skew-symmetric; this shows $\textup{(e)}\Rightarrow \textup{(f)}$.
Finally, the condition for $\ad_{\mf X}$ to be a skew-isometry for all $\mf X\in \frak g$ is
\begin{align}
    \langle \mf E_j, \ad_{\mf E_\ell}\mf E_i \rangle = -\langle \ad_{\mf E_\ell}\mf E_j,\mf E_i \rangle \ \Leftrightarrow\ C_{\ell i}^k \mf W_{kj} = - C_{\ell j}^k \mf W_{ki} 
\end{align}
for all $i,j,\ell \in \lbrace 1,\ldots,n\rbrace$, which shows $\textup{(f)}\Leftrightarrow\textup{(b)}$ and completes the proof.
\end{proof}

If one chooses to work with an orthonormal basis, then the conditions in \cref{lem:bi-invariant-metrics} take a simpler form, as presented in the following corollary.
\begin{corollary}
If the basis $\lbrace \mf E_i \rbrace_{i=1}^n$ is orthonormal (i.e., $\mf W=\mf I_n$), then conditions \textup{(a)} and \textup{(b)} of Lemma \ref{lem:bi-invariant-metrics} are equivalent to 
\begin{align}
\bm\Ad_{\mf g}\tran \bm\Ad_{\mf g} =\mf I_n\quad\text{and}\quad\bm\ad_\mf X = - \bm\ad_\mf X\tran
\label{eq:ad-condition}
\end{align}
for all $g \in G$ and $\mf X\in \mathfrak g$, respectively. That is, $\bm\Ad_{\mf g}$ is orthogonal and $\bm\ad_\mf X$ is skew-symmetric. Equivalently, the structure constants satisfy $C_{\ell i}^j = -C_{\ell j}^i$ for all $i,j,\ell \in \lbrace 1,\ldots,n\rbrace$.
\end{corollary}

\subsection{Geodesics}\label{app:geodesics}
Let $G$ be endowed with a left-invariant (but not necessarily bi-invariant) Riemannian metric. 
Consider a smooth curve $\bm\upgamma:[0, 1] \rightarrow G$ on $G$. Given $\mf g_1,\mf g_2\in G$, the curve $\bm\upgamma$ is the unique \textit{geodesic} connecting $\mf g_1$ and $\mf g_2$ if and only if it is a critical point of the \textit{energy functional}
\begin{align}
\mathcal E(\bm\upgamma)&\coloneqq \int_0^1 {\ip{\dot{\bm\upgamma}(t)}{\dot{\bm\upgamma}(t)}^L_{\bm\upgamma(t)}}dt,
\end{align} 
subject to the constraints on $\bm\upgamma(0)=\mf g_1$ and $\bm\upgamma(1)=\mf g_2$ \cite[p. 189]{lee2018introduction}.\footnote{
    While the minimization of $\mathcal E$ completely characterizes geodesics (and therefore, is sufficient for our purposes), the \textit{definition} of a geodesic is typically stated in terms of the Levi Civita connection of the Riemannian metric \cite{lee2018introduction}. The Levi Civita connection of a left-invariant metric is discussed in \cite{cheegerComparisonTheoremsRiemannian1975a} and \cite{park1995distance}.} Here, $\dot{\bm\upgamma}(t)\coloneqq \frac{d}{ds}\bm\upgamma(s)\big|_{s=t}$ represents the velocity of $\bm\upgamma$ at $t$.\footnote{Intrinsically, the velocity vector is given by a pushforward:
$\dot{\bm\upgamma}(t) \coloneqq d\bm\upgamma_t\left(\frac{\partial}{\partial s}\big|_{s=t}\right)$, where $\frac{\partial}{\partial s}$ is the unit vector field on $[0,1]$.}


Letting $\bm\upxi(t)\coloneqq \bm\upgamma(t)^{-1}\dot{\bm\upgamma}(t)$, we 
can use the \textit{Euler-Poincaré equations} \cite[Thm. 13.5.3]{marsden2013introduction} to show that $\bm\upgamma$ is a geodesic if and only if \cite{modin2010geodesics,guigui2023introduction}
\begin{align}
\frac{d}{dt}\bm\upxi^\vee(t)&= \mf W^{-1} \bm\ad_{\bm\upxi(t)}\tran \mf W \bm\upxi^\vee(t).
\end{align}
In \cite{guigui2023introduction}, the same equation is written as $\frac{d}{dt}\bm\upxi^\vee(t)=\bm\ad_{\bm\upxi(t)}^\ast \bm\upxi^\vee(t)$. As shown in the proof of \cref{lem:bi-invariant-metrics}, the matrix of $\ad_{\bm\upxi(t)}^\ast$ is given by $\bm \ad_{\bm\upxi(t)}^\ast = \mf W^{-1}\bm\ad_{\bm\upxi(t)}\tran \mf W$, so these are indeed equivalent expressions.
Unlike the geodesic equation on general Riemannian manifolds, which is a system of second-order differential equations, the variational approach (based on the Euler-Poincaré equations) yields first-order differential equations.

If the inner product on $\mathfrak g$ is $\Ad$-invariant, then Lemma \ref{lem:bi-invariant-metrics} says that $\bm\ad_{\bm\upxi(t)}\tran \mf W \bm\upxi^\vee(t)=\mf 0$. Consequently, 
$\frac{d}{dt}\bm\upxi^\vee(t) = \mf 0$, verifying that the geodesics of a bi-invariant metric are constant-velocity curves of the form $\bm\upgamma(t)=g\exp\big(t\bm\upxi(0)\big)$. If there does not exist an $\Ad$-invariant inner product on $\mathfrak g$, then the geodesic equation is solved by a time-varying curve of the form $\bm\upxi:[0,1]\rightarrow \mathfrak g$.


\bibliographystyle{siamplain}
\bibliography{references}
\end{document}


\maketitle

\section{A detailed example}

Here we include some equations and theorem-like environments to show
how these are labeled in a supplement and can be referenced from the
main text.
Consider the following equation:
\begin{equation}
  \label{eq:suppa}
  a^2 + b^2 = c^2.
\end{equation}
You can also reference equations such as \cref{eq:matrices,eq:bb} 
from the main article in this supplement.

\lipsum[100-101]

\begin{theorem}
An example theorem.
\end{theorem}

\lipsum[102]
 
\begin{lemma}
An example lemma.
\end{lemma}

\lipsum[103-105]

Here is an example citation: \cite{KoMa14}.

\section[Proof of Thm]{Proof of \cref{thm:bigthm}}
\label{sec:proof}

\lipsum[106-112]

\section{Additional experimental results}
\Cref{tab:smfoo} shows additional
supporting evidence. 

\begin{table}[htbp]
\footnotesize
  \caption{Example table.}\label{tab:smfoo}
\begin{center}
  \begin{tabular}{|c|c|c|} \hline
   Species & \bf Mean & \bf Std.~Dev. \\ \hline
    1 & 3.4 & 1.2 \\
    2 & 5.4 & 0.6 \\ \hline
  \end{tabular}
\end{center}
\end{table}

\bibliographystyle{siamplain}
\bibliography{references}